\documentclass{amsart}

\usepackage[mathcal]{euscript}
\usepackage{amssymb, amsfonts, xypic}

\newtheorem{theorem}{Theorem}[section]
\newtheorem{lemma}[theorem]{Lemma}
\newtheorem{corollary}[theorem]{Corollary}
\newtheorem{proposition}[theorem]{Proposition}

\theoremstyle{definition}
\newtheorem{definition}[theorem]{Definition}

\theoremstyle{remark}
\newtheorem{remark}[theorem]{Remark}

\numberwithin{equation}{section}

\newcommand{\Z}{\mathbb{Z}}
\newcommand{\N}{\mathbb{N}}

\newcommand{\C}{\mathcal{C}}
\newcommand{\LS}{\mathcal{LS}}
\newcommand{\SSet}{\mathcal{SS}}
\newcommand{\SFS}{\mathcal{SFS}}

\DeclareMathOperator{\Ho}{Ho}
\DeclareMathOperator{\Hom}{Hom}

\DeclareMathOperator*{\colim}{colim}
\DeclareMathOperator*{\holim}{holim}

\newcommand{\Ab}{\textup{Ab}}
\newcommand{\strict}{\textup{strict}}
\newcommand{\pro}{\textup{pro-}}
\newcommand{\pros}{\textup{pro}}
\newcommand{\Map}{\textup{Map}}
\newcommand{\fin}{\textup{fin}}

\newcommand{\ol}{\overline}
\newcommand{\map}{\rightarrow}

\begin{document}

\title[A model structure for pro-simplicial sets]
{A model structure on the category of pro-simplicial sets}

\author{Daniel C.\,Isaksen}
\address{
Fakult\"at f\"ur Mathematik\\
Universit\"at Bielefeld\\
33501 Bielefeld, Germany}
\email{isaksen@mathematik.uni-bielefeld.de}
\thanks{The author was supported in part by an NSF Graduate Fellowship.}

\subjclass{Primary 18E35, 55Pxx, 55U35; Secondary 14F35, 55P60}
\date{June 18, 2001}
\keywords{Closed model structures, pro-spaces, \'etale homotopy}

\begin{abstract}
We study the category $\pro \SSet$ of pro-simplicial sets, which arises
in \'etale homotopy theory, shape theory,
and pro-finite completion.  We establish a
model structure on $\pro \SSet$ so that it is possible to do homotopy
theory in this category.  This model structure is closely related to the
strict structure of Edwards and Hastings.  
In order to understand the notion of homotopy
groups for pro-spaces we use local systems on pro-spaces.
We also give several alternative descriptions of weak equivalences, including
a cohomological characterization.  We outline dual constructions for
ind-spaces.
\end{abstract}

\maketitle

\section{Introduction}
\label{sctn:intro}

If $\C$ is a category, then the pro-category $\pro\C$
\cite[Expos\'e\, 1, Section 8]{SGA} is
the category whose objects are small cofiltered systems in $\C$ (of arbitrary
shape) and whose morphisms are given by the formula

$$\Hom(X, Y) = \lim_s \colim_t \Hom_{\C}(X_t, Y_s).$$ While investigating
the \'etale homotopy functor, Artin and Mazur \cite{AM} studied the
category $\pro \Ho(\SSet)$, where Ho$(\SSet)$ is the homotopy category of
simplicial sets.  They also introduced a notion of weak
equivalence of pointed connected pro-spaces that involved isomorphisms of
pro-homotopy groups.

However, an Artin-Mazur weak equivalence is not the same
as an isomorphism in $\pro \Ho(\SSet)$.
This suggests that $\pro \Ho(\SSet)$ is
not quite the correct category for studying \'etale homotopy.

Around the same time, Quillen \cite{DQ} developed the fundamental notions
of homotopical algebra by realizing that model structures 
allow one to do homotopy theory in many different categorical contexts.
A model structure on a category is a choice of three classes of maps
(weak equivalences, cofibrations, and fibrations) satisfying certain
axioms.  The weak equivalences are inverted to obtain the associated
homotopy category, while
the cofibrations and fibrations serve an auxiliary role.
Quillen \cite[II.0.2]{DQ} observed that the homotopy theory
of pro-spaces would be an interesting application of model structures.

At least two model structures for pro-spaces are already known to exist.
Edwards and
Hastings \cite{EH} established a ``strict'' model structure on pro-spaces for
the purposes of shape theory and proper homotopy theory, but their
weak equivalences did not generalize those of Artin and Mazur.

Also, 
Grossman \cite{JG} described a different
model structure for pro-spaces that are countable towers.  The weak
equivalences of this theory are appropriately related to the Artin-Mazur
equivalences.
The category of towers is suitable in many applications
from proper homotopy theory
because it is reasonable to assume that the neighborhoods
at infinity have a countable basis.  

However, 
applications of pro-spaces to the algebro-geometric concept
of \'etale cohomology require
more general pro-spaces.  General cofiltered
systems of spaces are necessary for essentially the same reason that
the sheaf theory of Grothendieck topologies is necessary 
to define \'etale cohomology.

This paper gives a generalization of Grossman's model structure
to the entire category of pro-spaces.
Weak equivalences between pointed connected pro-spaces
are precisely Artin-Mazur weak equivalences.  
The Edwards-Hastings strict structure is an intermediate stage to
building our structure.
Our homotopy theory is the $P$-localization of the strict homotopy theory,
where $P$ is the functor that replaces a space with its Postnikov tower.

Our weak equivalences have
several alternative characterizations, one of which uses
twisted cohomology.  This
is important because it is often difficult to check that a map of pro-spaces
induces an isomorphism on homotopy groups.  Usually it is much easier to
verify a cohomology isomorphism and then conclude that the map is a weak
equivalence.

Another characterization of weak equivalences is in terms of ``eventually
$n$-connected'' maps (see Theorem \ref{thm:tfae} {\em (d)}), which 
is a very convenient property in practice.
The equivalence of this property with the definition of weak equivalence
is not obvious.  One must use the full power of the model structure to prove
the equivalence.

We mention two applications of this homotopy theory of pro-spaces.
First, the model structure gives a more
convenient category for studying \'etale homotopy because it allows
the reinterpretation of the central ideas of the theory in terms of the
established notions of model structures. 
It is also a start towards the definition of generalized cohomology
of pro-spaces.  For example, define $K^0$
to be represented by the constant pro-space $BU$.  The realization of
$K$-theory as a generalized cohomology theory requires more understanding
of the category of pro-spectra.
The \'etale 
$K$-theory of a scheme \cite{etK} is probably most clearly expressed
as a generalized cohomology theory applied to the \'etale homotopy
type.

Second, pro-spaces arise in the study of pro-finite completion
\cite{DS} \cite{BK} \cite{FM}.  Again, the new model structure provides a
better category in which to study such completions.  For example, Mandell
\cite{MAM} has used the model structure to compare his new algebraic 
construction of
$p$-adic homotopy theory to Sullivan's $p$-pro-finite completion.

We describe briefly the model structure; the formal definitions appear in
Section \ref{sctn:structure}.
Given a pointed pro-space $X$ ({\em i.e.,} a map from the one-point
constant pro-space $*$ to $X$ or equivalently a pro-object in the category
of pointed spaces), one can apply the functor $\pi_n(-,*)$ to
each space in the pro-system to get a pro-group $\pi_n(X,*)$. 
The most obvious notion of equivalence of pro-spaces $X \map Y$ is the
requirement that the map induce an isomorphism of pro-homotopy groups
$\pi_n(X,*) \map \pi_n(Y,*)$ for every $n \geq 0$ and every point $*$ of
$X$. 

However, points of a pro-space are rather awkward.  In fact, some non-trivial
pro-spaces have {\em no} points whatsoever.  
Local systems permit the discussion of
homotopy groups without choosing basepoints.
Grossman's work on towers \cite{JG} inspired this trick.

Cofibrations of pro-spaces are maps that are isomorphic to levelwise
cofibrations of systems of spaces of the same shape.  The model category
axioms then force the definition of fibrations.
These fibrations are similar to those of 
Edwards and Hastings \cite[3.3]{EH}, but they satisfy an extra
condition that compares the homotopy groups of the total pro-space to the 
homotopy groups of the base pro-space.  

The model structure has a few interesting aspects.  For example,
it is not cofibrantly generated.
This means that the proof of the model axioms is quite different from
the standard arguments.  One result of this fact is that the
factorizations are not functorial.  Almost all naturally arising model
structures are cofibrantly generated.  Hence pro-spaces are an interesting
example of a non-cofibrantly generated model structure.

Also, we only
know how to work with pro-simplicial sets, not pro-topological spaces.
We use in a very significant way the fact that finite dimensional skeletons
are functorial
for simplicial sets.  Since relative cell complexes do not have functorial
skeletons, the same proofs do not apply.

Finally, note that the
category of ind-spaces has a similar model structure.  We do not provide
details in this paper because we have no application in mind.  Several
comments throughout the paper explain where the significant differences
occur.

The paper is divided into three main parts.  Sections 
\ref{sctn:prelim-pro}--\ref{sctn:htpy}
give some background material and introduce language and tools
for later use.  Sections \ref{sctn:structure}--\ref{sctn:strict}
describe the model structure,
state the main theorems, and make comparisons to other homotopy theories.
Sections \ref{sctn:colim}--\ref{sctn:cofibgen}
provide proofs of the main theorems.  

We assume familiarity with model structures.  See \cite{PH} or
\cite{DQ} for background material.

Many of the important ideas in this paper come from Grossman \cite{JG}.
I thank Peter May, Bill Dwyer, Brooke Shipley, Greg Arone, Michael Mandell,
and Charles Rezk for many
helpful conversations throughout the progress of this work.

\section{Preliminaries on Pro-Categories} 
\label{sctn:prelim-pro}

First we establish some terminology for pro-categories.

\begin{definition} 
\label{defn:pro}
For a category $\C$, the category $\pro \C$ has objects all small cofiltering
systems in $\C$, and 
$$\Hom_{\pro \C}(X,Y) = \lim_s \colim_t \Hom_{\C}
     (X_t, Y_s).$$
\end{definition}

Objects of $\pro \C$ can be thought of as functors from arbitrary small
cofiltering categories to $\C$.  See \cite{SGA}
or \cite{AM} for more background on the definition of pro-categories.
We use both set theoretic and categorical
language to discuss indexing categories;
hence ``$t \geq s$'' and ``$t \map s$'' mean the same thing.

The word ``system'' always refers to an object of a pro-category,
while the word ``diagram'' refers to a diagram of pro-objects.

A {\em subsystem}
of an object $X: I \map \C$ in $\pro \C$ is a restriction of $X$
to a cofiltering subcategory $J$ of $I$.  A subsystem is {\em cofinal} if for
every $s$ in $I$, there exists some $t$ in $J$ and an arrow $t \map s$ in $I$.
A system is isomorphic in $\pro \C$ to any of its cofinal subsystems.

A directed set $I$ is {\em cofinite}
if for every $t$, the set of elements of $I$ less than $t$ is finite.  
Except in Section \ref{sctn:colim},
all systems are indexed by
cofinite directed sets rather than arbitrary cofiltering categories.
This is no loss of generality \cite[2.1.6]{EH}
(or \cite[Expos\'e\, 1, 8.1.6]{SGA}).  
The cofiniteness is critical 
because many constructions and proofs proceed inductively.

Whenever possible we avoid
mentioning the structure maps of a pro-object $X$ explicitly.
When necessary the notation $(X, \phi)$ indicates a system of
objects $\{ X_s \}$ with structure maps $\phi_{ts}: X_t \map X_s$.

\begin{definition} 
\label{defn:height}
Let $I$ be a cofinite directed set.  For each $s$ in $I$,
the {\em height} $h(s)$ of $s$ is the value of $n$ in the longest chain
$s > s_1 > s_2 > \cdots > s_n$ starting at $s$ in $S$.  In particular,
$h(s) = 0$ if and only if there are no elements of $I$ less than $s$.
\end{definition}

We always assume that systems have no initial object or equivalently that
any system has objects of arbitrarily large height.  This is no loss
of generality
since it is possible to add isomorphisms to a system so that it no
longer has an initial object, and this new system is isomorphic to the old
one.

We frequently consider maps between two pro-objects with the same
index categories.  In this setting, a {\em level map}
$X \map Y$ between pro-objects
indexed by $I$ is given by maps $X_s \map Y_s$ for all $s$ in $I$.
Up to isomorphism, every map is a level map \cite[Appendix 3.2]{AM}.

A map satisfies a certain property {\em levelwise} if it is a level map
$X \map Y$ such that each $X_s \map Y_s$ satisfies that property.

\begin{lemma} 
\label{lem:pro-iso}
A level
map $A \map B$ in $\pro \C$ is an isomorphism if and only if for all $s$,
there exists $t \geq s$ and a commutative diagram

$$\xymatrix{
A_t \ar[r] \ar[d] & B_t \ar[d]\ar[dl] \\
A_s \ar[r]        & B_s.              }$$
\end{lemma}

\begin{proof}
The maps $B_t \map A_s$ induce an inverse.
\end{proof}

\section{Preliminaries on Simplicial Sets} 
\label{sctn:prelim-ss}

Now we review some definitions and results about simplicial sets.  Let
$\SSet$ be the category of simplicial sets.  We use the expressions
``space'' and ``simplicial set'' interchangeably.

For simplification, we often use the
same notation for a basepoint and its image under various maps ({\em e.g.},
$\pi_n(X,*) \map \pi_n(Y,*)$).

\begin{definition} 
\label{defn:equiv}
A map $f: X \map Y$
of simplicial sets is an {\em $n$-equivalence} if for all 
basepoints $*$ in $X$, $f$ induces an isomorphism
$\pi_i(X,*) \map \pi_i(Y,*)$ for $0 \leq i < n$ and
a surjection $\pi_n(X,*) \map \pi_n(Y,*)$.
The map $f$ is a {\em co-$n$-equivalence} if
for all basepoints $*$ in $X$, $f$ induces an isomorphism
$\pi_i(X,*) \map \pi_i(Y,*)$
for $i > n$ and an injection $\pi_n(X,*) \map \pi_n(Y,*)$.
\end{definition}

\begin{definition} 
\label{defn:n}
Set $$J_n = \{ \Lambda^m_k \map \Delta^m | m \geq 0 \} \cup
           \{ \partial \Delta^m \map \Delta^m | m > n \}.$$
A map of simplicial sets is a {\em co-$n$-fibration}
if it has the right lifting
property with respect to all maps in $J_n$.
A map of simplicial sets is an {\em $n$-cofibration}
if it has the left lifting property with respect to all co-$n$-fibrations.
\end{definition}

In other words, a co-$n$-fibration is a $J_n$-injective, and an
$n$-cofibration
is a $J_n$-cofibration \cite[12.4.1]{PH}.  Note that co-$n$-fibrations and
$n$-cofibrations are characterized by lifting properties with respect to
each other.
When $n = -1$ or $n = \infty$,
the definitions reduce to the usual definitions of trivial cofibrations
and fibrations or to the definitions of cofibrations and trivial fibrations.

We will see below that $n$-cofibrations are just maps that are
both cofibrations and $n$-equivalences.  Also, co-$n$-fibrations are just
maps that are both fibrations and co-$n$-equivalences.  
These facts motivate the terminology.

\begin{proposition} 
\label{prop:SS-factor}
For any $n$, each map $f$ of simplicial sets factors as $f = pi$, where
$i$ is an $n$-cofibration and $p$ is a co-$n$-fibration.
\end{proposition}

\begin{proof}
Apply the small object argument \cite[12.4]{PH}.
\end{proof}

\begin{lemma} 
\label{lem:co-n-fib}
A map $f:E \map B$
is a co-$n$-fibration if and only if $f$ is a fibration and
co-$n$-equivalence.
\end{lemma}

\begin{proof}

First suppose that $f$ is a co-$n$-fibration. 
The generating trivial cofibrations are contained in $J_n$, so $f$ is
also a fibration.  Now we must show that $f$ is also a co-$n$-equivalence.

Consider test diagrams of the form
$$\xymatrix{
\partial \Delta^m \ar[r]\ar[d] & E  \ar[d]^f \\
\Delta^m \ar[r]_g \ar@{-->}[ur] & B,    }$$
where $m > n$.
Let $g:\Delta^m \map B$ be a constant map with image $*$, and let $F$
be the fiber of $f$ over $*$.
Since lifts exist in the test diagram,
$\pi_m F = 0$ for $m \geq n$.  
Using the long exact sequence of homotopy groups of a
fibration, it follows that $f$ is a co-$n$-equivalence.

Now suppose that $f$ is a fibration and co-$n$-equivalence.  It follows
from the long exact sequence of homotopy groups
that $\pi_m F = 0$ for $m \geq n$, where $F$ is any fiber of $f$.
There are lifts in the test diagrams shown above
for $m > n$ because the obstructions to such lifts belong to $\pi_{m-1} F$.
Hence $f$ is a co-$n$-fibration.
\end{proof}

\begin{definition} 
\label{defn:relative-sk}
If $A \map X$ is a cofibration of simplicial sets, then the
{\em relative $n$-skeleton} $X^{(n)}$ is the union of $A$ and the
$n$-skeleton of $X$.
\end{definition}

\begin{lemma} 
\label{lem:n-cofib}
A map $f:A \map X$ is an $n$-cofibration if and only if
$f$ is a cofibration and $n$-equivalence.
\end{lemma}

\begin{proof}
Since trivial fibrations are co-$n$-fibrations, 
$n$-cofibrations are also cofibrations.

Suppose that $f: A \map X$ is a relative
$J_n$-cell complex \cite[12.4.6]{PH}.  
Then $A \map X^{(n)}$ is a weak
equivalence since $X^{(n)}$ is obtained from $A$
by gluing along maps of the form $\Lambda_k^m \map \Delta^m$.
Note that $X^{(n)} \map X$ is an $n$-equivalence, so $A \map X$ is
also an $n$-equivalence.  Arbitrary $n$-cofibrations are retracts
of relative $J_n$-cell complexes \cite[13.2.9]{PH}, so all
$n$-cofibrations are $n$-equivalences.

Conversely, suppose that $f$ is a cofibration and $n$-equivalence.
We show that $f$ is an $n$-cofibration by demonstrating that it
has the left lifting property with respect to all maps that are both
fibrations and co-$n$-equivalences.  By Lemma \ref{lem:co-n-fib}, this
means that $f$ has the left lifting property with respect to all
co-$n$-fibrations.

Factor $f$ as

$$\xymatrix@1{A \ar[r]^j & Y \ar[r]^q & X,    }$$
where $j$ is a trivial cofibration and $q$ is a fibration.
Note that $q$ is also an $n$-equivalence.
Let $p: E \map B$ be a map that is both a fibration and co-$n$-equivalence,
and consider a diagram

$$\xymatrix{
A \ar[rr]\ar[d]_j & & E \ar[d]^p \\
Y \ar@{-->}[urr] \ar[r]_q & X \ar[r] & B.  }$$
The dashed arrow exists because $j$ is a trivial cofibration and
$p$ is a fibration.  

This gives the diagram

$$\xymatrix{
A \ar[r]\ar[d]_f & Y \ar[r] & E \ar[d]^p \\
X \ar[r]       & X \ar[r] &  B.    }$$

There is no obstruction to lifting over $X^{(0)}$ since
$\pi_0 A \map \pi_0 X$ is surjective.

Obstructions to finding a lift over
the higher relative skeletons of $X$ belong to
$\pi_m F$, where $F$ is some fiber of $p$.  These obstructions lie
in the image of $\pi_m G$, where $G$ is a fiber of $q$.  For every $m$,
either $\pi_m G$ or $\pi_m F$ is zero.  Hence there are no obstructions
to lifting.
\end{proof}

\begin{remark}
For ind-spaces, consider the set
$$I_n = \{ \Lambda^m_k \map \Delta^m | m \geq 0 \} \cup
           \{ \partial \Delta^m \map \Delta^m | m < n \}$$
to define {\em $n$-fibrations} and {\em co-$n$-cofibrations}.
All $n$-fibrations are both fibrations
and $n$-equivalences, but the converse is not true.
If $f: A \map X$ is a co-$n$-cofibration, then $f$ is a cofibration
and the induced map $X^{(n-1)} \map X$ is a weak equivalence.
\end{remark}

\section{Local Systems} 
\label{sctn:loc-sys}

The language of local systems is necessary in order to state the idea
of homotopy groups for pro-spaces.  Recall that a local system on 
a space $X$ is a functor $\Pi X \map \Ab$, where $\Pi X$ is the fundamental
groupoid of $X$ and $\Ab$ is the category of abelian groups.  
Denote by $\LS (X)$ the
category of local systems on $X$ or equivalently the category of
locally constant sheaves on $X$.

For example, $\Pi_nX$
is a local system on $X$ for $n \geq 2$.  It is defined
by $(\Pi_nX)_x = \pi_n(X,x)$
with isomorphisms $(\Pi_nX)_x \map (\Pi_nX)_y$ given by the usual maps
on homotopy groups induced by paths.  

Occasionally we refer to local systems with values in non-abelian groups.
For example, $\Pi_1 X$ is such a local system.  
We emphasize the notational 
distinction between $\Pi X$ (a groupoid) and $\Pi_1 X$ (a local system).

If $f:X \map Y$ is a map of spaces, then $f$ induces a map of local systems 
$\Pi_nX \map f^*\Pi_nY$ on $X$, where $f^*$ is the pullback functor
$\LS(Y) \map \LS(X)$.
Recall that the functor $f^*$ is exact in the sense that it preserves finite
limits and colimits, and it is also exact in the sense that
it preserves exact sequences.  

\begin{lemma} 
\label{lem:SS-we}
Let $f: X \map Y$ be a map of spaces.
Then $f$ is a weak equivalence if and only if 
$\pi_0 f$ is an isomorphism and $\Pi_nX \map f^*\Pi_nY$ is an
isomorphism of local systems on $X$ for all $n \geq 1$.
\end{lemma}

\begin{proof}
This is a restatement without reference to basepoints of the
definition of weak equivalence of simplicial sets.
\end{proof} 

We now extend the definitions to pro-spaces.

\begin{definition} 
\label{defn:pro-loc-sys}
A {\em local system} on a pro-space $X$ is an object of
$\colim_s \LS (X_s)$.  If $L$ is a local system on $(X, \phi)$
represented by a functor $L_s: \Pi X_s \map \Ab$
and $M$ is another local system on $X$ represented by a functor
$M_t: \Pi X_t \map \Ab$, then a map from $L$ to $M$ is an
element of $\colim_u \Hom_{\LS (X_u)} (\phi^*_{us} L_s, \phi^*_{ut} M_t)$.
Denote by $\LS (X)$ the category of local systems on $X$.
\end{definition}

A local system on $X$ is represented by a local system on $X_s$
for some $s$.  For example, for $n \geq 1$, $\Pi_nX_s$
is a local system on $X$ for each $s$.  

A map between two local systems on $X$ is
a map of representing local systems pulled back to some $X_u$.
For example, for $n \geq 1$ and $t \geq s$,  
$\Pi_n X_t \map \Pi_n X_s$ is a map of local systems on $(X,\phi)$
because $\Pi_n X_t \map \phi^*_{ts} \Pi_n X_s$ is a map of local systems
on $X_t$.

Let $f:X \map Y$ be a map of pro-spaces, and let
$L: \Pi Y_s \map \Ab$ be a local system on $Y_s$.  Choose any
map $f_{ts}: X_t \map Y_s$ representing
$f$ and consider the
functor $L \circ f_{ts}: \Pi X_t \map \Pi Y_s \map \Ab$.  This gives
a well-defined functor $\colim_s \LS(Y_s) \map \colim_t \LS(X_t)$.

\begin{definition} 
\label{defn:pro-loc-sys-pullback}
Let $f:X \map Y$ be a map of pro-spaces.  The {\em pullback} 
$f^*: \LS(Y) \map \LS(X)$ is the functor
$\colim_s \LS(Y_s) \map \colim_t \LS(X_t)$.
\end{definition}

\begin{lemma} 
\label{lem:pro-pullback-exact}
If $f: X \map Y$ is a
map of pro-spaces, then $f^*: \LS (Y) \map \LS (X)$ is an exact functor
in the sense that it preserves finite limits and colimits.
\end{lemma}

\begin{proof}
Without loss of generality, we may assume that $f$ is a level map.
Given a finite diagram of local systems $L$ on $Y$, there
is an $s$ such that each $L^i$ is represented by a local system
$L^i_s$ on $Y_s$.  
Then for each $i$, $f^*L^i$ is represented by $f_s^*L^i_s$.
Now $\colim L^i$ in $\LS (Y)$ is represented by $\colim L_s^i$ in
$\LS (Y_s)$, so
$f^* (\colim L^i)$ is represented by $f_s^* (\colim L^i_s)$.
Also, $\colim f^*L^i$ in $\LS (X)$ is represented by
$\colim f_s^*L_s^i$ in $\LS (X_s)$.
But $f_s^*$ commutes with finite colimits, so $\colim f^* L^i$ and
$f^* (\colim L^i)$ are isomorphic since they are represented by the same
local system on $X_s$.

An identical argument shows that $f^*$ commutes with finite limits.
\end{proof}

It follows from the lemma that $f^*$ is exact in the sense that
it preserves exact sequences.

\section{Homotopy Groups} 
\label{sctn:htpy}

With the notions of local systems in place, 
we can define homotopy groups of pro-spaces as pro-objects in a
category of local systems.
The local systems are necessary to avoid mentioning basepoints.

\begin{definition} 
\label{defn:pro-htpy-gp}
If $X$ is a pro-space and $n \geq 2$,
then $\Pi_nX$ is the pro-local system on $X$ given
by $\{ \Pi_nX_s \}$.  Also, $\pi_0 X$ is the pro-set given by
$\{ \pi_0 X_s \}$, and $\Pi_1 X$ is the pro-local system on $X$ with
values in non-abelian groups given by $\{\Pi_1 X_s \}$.
\end{definition}

Note that a map of pro-spaces $f: X \map Y$ 
induces a map $\Pi_nX \map f^*\Pi_nY$ in $\pro \LS (X)$.

\begin{lemma} 
\label{lem:pro-pullback-pro-iso}
If $f: X \map Y$ is a map of pro-spaces and $\pi_0 f$ is an epimorphism in
the category of pro-sets, then
a map $g$ of pro-local systems is an isomorphism in $\pro \LS(Y)$
if and only if $f^*(g)$ is an isomorphism in $\pro \LS(X)$.
\end{lemma}

\begin{proof}
Without loss of generality, assume that $f$ is a level map.
Note that $\pro \LS(Y)$ is an abelian category \cite[Appendix 4.5]{AM}.
Since $f^*$ is exact, it suffices to 
consider a local system $L$ on $Y$ such that $f^* L = 0$ and conclude that
$L = 0$.  

A pro-local system $M$ is zero if and only if for every $i$, there
exists $j \geq i$ such that the map $M^j \map M^i$ is trivial.
For any $i$, choose $j \geq i$ so that $f^*L^j \map f^*L^i$ is trivial.
Let $L^i_s$ and $L^j_s$ be local systems
on $Y_s$ representing $L^i$ and $L^j$ respectively.
Choose $s$ large enough so that $f^*_s L^j_s \map f^*_s L^i_s$ is a trivial
map of local systems on $X_s$.

Since $\pi_0 f$ is an epimorphism, there exists some $t$ and a map 
$\pi_0 Y_t \map \pi_0 X_s$ such that the map $\pi_0 Y_t \map \pi_0 Y_s$
factors through $\pi_0 X_s$.
Since $f_s^* L^j_s \map f_s^* L^i_s$ is trivial, the map
$L^j_s \map L^i_s$ is trivial when
restricted to the components of $Y_s$ in the image of $X_s$.  Now
the image of $\pi_0 X_s$ in $\pi_0 Y_s$ contains the image of $\pi_0 Y_t$,
so $L^j_s \map L^i_s$ becomes trivial when pulled back to $Y_t$.
Hence the map $L^j \map L^i$ is trivial.  This means that the pro-local
system $L$ is zero.
\end{proof}

\begin{remark}
A similar statement is true for pro-local systems with values in
non-abelian groups.
\end{remark}

\begin{lemma} 
\label{lem:levelwise-components}
A map of pro-spaces $f: X \map Y$ induces an isomorphism of pro-sets
$\pi_0 f: \pi_0 X \map \pi_0 Y$ if and only if $f$ is isomorphic to a level map
$f': X' \map Y'$ such that $f'$ induces a level isomorphism
$\pi_0 f': \pi_0 X' \map \pi_0 Y'$.
\end{lemma}

\begin{proof}
One direction is clear because level isomorphisms are pro-isomorphisms.
Assume that $\pi_0 f$ is an isomorphism.  We may
also assume that $f$ is a level map.

Define $X' = X$ and $Y' = Y \times_{\pi_0 Y} \pi_0 X$.  Here
we identify pro-sets with pro-spaces of dimension zero.

Then
$Y'$ is isomorphic to $Y$ since $\pi_0 X \map \pi_0 Y$ is an isomorphism.
Let $f': X' \map Y'$ be the map induced by $f: X \map Y$ and 
the projection $X \map \pi_0 X$.
Pullbacks can be constructed levelwise in pro-categories,
so for all $s$, 
$Y'_s = Y_s \times_{\pi_0 Y_s} \pi_0 X_s$ and $f'_s$ is induced by 
$f_s: X_s \map Y_s$ and $X_s \map \pi_0 X_s$.
Note that $f'_s$ induces an isomorphism $\pi_0 X'_s \map \pi_0 Y'_s$.
Hence $f'$ is the desired map.
\end{proof}

If $f: X \map Y$ is a map of spaces such that $\pi_0 f$ is an isomorphism
and $\pi_1 f$ is an isomorphism for every basepoint, then $f^*$
induces an equivalence of categories $\LS(Y) \map \LS(X)$.  The following
lemma makes an analogous statement for pro-spaces.

\begin{lemma} 
\label{lem:pullback-surj}
If $f: X \map Y$ is a map of pro-spaces such that
$\pi_0 f$ and $\Pi_1 X \map f^* \Pi_1 Y$ are isomorphisms,
then the functor $f^*: \LS (Y) \map \LS (X)$ is an equivalence of 
categories.
\end{lemma}

\begin{proof}
We only prove that $f^*$ is essentially surjective in the sense that
every object $L$ of $\LS (X)$ is isomorphic to $f^* M$ for 
some $M$ in $\LS (Y)$.
We leave the rest of the proof to the interested reader.  We will use only
the surjectivity in this work.

With no loss of generality, we may assume that $f$ is a level map.  By Lemma
\ref{lem:levelwise-components}, we may also assume that
$\pi_0 f$ is a level isomorphism.

Let $L$ be a local system on $(X, \phi)$ represented by
a local system $L_s$ on $X_s$.  There
exists $t \geq s$ and a commutative diagram

$$\xymatrix{
\Pi_1 X_t \ar[r]\ar[d] & f_t^*\Pi_1 Y_t \ar[d]\ar[dl] \\
\phi^*_{ts} \Pi_1 X_s \ar[r] & \phi^*_{ts} f_s^* \Pi_1 Y_s  }$$
of local systems on $X_t$.

Choose one point $x_i$ in each component of $X_t$.  Let
$y_i$ be the image of $x_i$ in $Y_t$; this is a choice of one point
in each component of $Y_t$.
Let $x'_i$ be the image of $x_i$ in $X_s$.

By evaluating the above diagram at $x_i$, the map
$\pi_1(X_t,x_i) \map \pi_1(X_s,x'_i)$ factors as
$$\xymatrix@1{
\pi_1(X_t,x_i) \ar[r] & \pi_1(Y_t,y_i) \ar[r]^{g_i} & \pi_1(X_s,x'_i).  }$$
The maps $g_i:\pi_1(Y_t, y_i) \map \pi_1(X_s,x'_i)$ and the local system $L_s$
determine a local system $M_t$ on $Y_t$ by setting 
$(M_t)_{y_i} = L_{x'_i}$ with
the $\pi_1(Y_t, y_i)$-action induced by $g_i$.  

Let $M$ be the
local system on $Y$ represented by $M_t$.  
Note that $f_t^* M_t$ is isomorphic to $\phi_{ts}^* L_s$.  Hence $f^* M$
is isomorphic to $L$.
\end{proof}

\begin{remark}
A similar statement applies to local systems with values in
non-abelian groups.
\end{remark}

\section{Model Structure} 
\label{sctn:structure}

Now we explicitly describe the model structure on the category of
pro-spaces.

\begin{definition} 
\label{defn:we} 
A map of pro-spaces $f:X \map Y$ is a {\em weak
equivalence} if $\pi_0 f$ is an isomorphism of pro-sets
and $\Pi_nX \map f^*\Pi_nY$ is an isomorphism in $\pro \LS (X)$
for all $n \geq 1$. 
\end{definition}

In Corollary \ref{cor:pointed-we} we will see that for
pointed connected pro-spaces,
a level map $X \map Y$ is a weak equivalence if and only if
$\pi_n X \map \pi_n Y$ is a pro-isomorphism for all $n$.
This works because there
is no need for arbitrary basepoints.  Hence Artin-Mazur weak
equivalences \cite[Section 4]{AM} are also weak equivalences.

\begin{definition} 
\label{defn:cofib}
A map of pro-spaces is a {\em cofibration}
if it is isomorphic to a levelwise cofibration.
\end{definition}

\begin{definition} 
\label{defn:fib}
A map of pro-spaces is a {\em fibration} if it has the right lifting
property with respect to all trivial cofibrations.
\end{definition}

\begin{theorem} 
\label{thm}
The above definitions give a proper simplicial
model structure on $\pro \SSet$ (without functorial factorizations).  This
model structure is not cofibrantly generated.
\end{theorem}

\begin{proof}
The axioms for a proper simplicial model structure are verified in
Sections \ref{sctn:colim} through \ref{sctn:proper}.

Limits and colimits exist by Proposition \ref{prop:complete}.  The
two-out-of-three axiom is Proposition \ref{prop:2/3}.
Retracts preserve weak equivalences because weak equivalences are defined
in terms of isomorphisms.  Retracts preserve fibrations because retracts
preserve lifting properties.
Corollary \ref{cor:cofib-retract} is the retract axiom for cofibrations.
Propositions \ref{prop:factor-trcofib-fib} and \ref{prop:factor-cofib-trfib}
are the factoring axioms, while Proposition \ref{prop:lift-cofib-trfib}
is the non-trivial lifting axiom.  The axioms for
a simplicial model structure are demonstrated in Proposition \ref{prop:simp}.
Proposition \ref{prop:proper}
shows that the model structure is proper, and Corollary \ref{cor:uncofibgen}
states that it is not cofibrantly generated.
\end{proof}

The model structure can be considered in two stages.  The
``strict'' structure of Edwards and Hastings \cite[3.3]{EH}, in which the
weak equivalences are defined levelwise, is an intermediate step;
see Section
\ref{sctn:strict} for details.  The situation is not unlike
the Bousfield-Friedlander 
strict and stable model structures for spectra \cite{BF}.

We assume Theorem \ref{thm} for the rest of this section and for the next
three sections.  Sections \ref{sctn:colim}--\ref{sctn:cofibgen} contain
the details of the proof of the theorem.

In practice we need a more concrete description of fibrations.
The next definition and proposition provide such a description.

\begin{definition} 
\label{defn:strong-fib}
A map is a {\em strong fibration} if it is isomorphic to a
level map of pro-spaces $X \map Y$ indexed by a cofinite
directed set such that for all $t$,
$$X_t \map Y_t \times_{\lim\limits_{s<t} Y_s} \lim\limits_{s<t} X_s$$
is a co-$n(t)$-fibration, where $n$ is any function from the index set
of $X$ to $\N$.
\end{definition}

Proposition \ref{prop:fib} shows that we 
could have taken the equivalent definition that 
$X_t \map Y_t \times_{\lim\limits_{s<t} Y_s} \lim\limits_{s<t} X_s$
is a fibration and $X_t \map Y_t$ is a co-$n(t)$-equivalence.  Formally
this alternative definition is more awkward, but 
it is often useful to recognize that specific examples
are strong fibrations.  

Note that this definition only applies to 
level maps indexed by cofinite directed sets.
In fact, the obvious generalization to systems of arbitrary shape 
is incorrect.  This suggests that
the above definition has a more natural equivalent formulation, but we do
not know how to restate it.

With no loss of generality,
we always assume that $n(t) \geq n(s)$ whenever $t \geq s$.

It is not obvious from the definition that strong fibrations are also
fibrations.  The technical heart of the whole theory lies in showing that
strong fibrations have the right lifting property with respect to 
trivial cofibrations.

\begin{proposition} 
\label{prop:fib-retract}
A map $f: E \map B$ of pro-spaces is a fibration if and only if
there exists a strong fibration
$q: E' \map B$ and a commutative diagram

$$\xymatrix{
E \ar[r] \ar[dr]_f & E' \ar[r] \ar[d]_q & E \ar[dl]^f \\
  & B              }$$
such that the composition along the top row is the identity on $E$.
\end{proposition}

\begin{remark}
The idea is that $f$ is a retract of $q$ in a strong sense where
the targets of $f$ and $q$ are the same.
\end{remark}

\begin{proof}
The ``if'' direction follows from the fact that strong fibrations are
fibrations by Proposition \ref{prop:lift-trcofib-strong-fib}
and the fact that retracts preserve fibrations.

Conversely, suppose that $f$ is a fibration.
Use Proposition \ref{prop:factor-trcofib-fib} to factor $f$ as

$$\xymatrix{ E \ar[r]^j & E' \ar[r]^q & B,}$$
where $j$ is a trivial cofibration and $q$ is a strong fibration.
The retract argument \cite[8.2.2]{PH} finishes the proof.
\end{proof}

\begin{proposition} 
\label{prop:trcofib-retract}
A level map 
$f: A \map X$ of pro-spaces is a trivial cofibration if and only if
there exists a levelwise cofibration
$i: A \map X'$ indexed by a cofinite directed set $I$
for which there is a strictly increasing function
$n:I \map \N$ such that $i_s$ is an $n(s)$-equivalence and there is
a commutative diagram

$$\xymatrix{
 & A \ar[dl]_f \ar[d]_i \ar[dr]^f \\
X \ar[r] & X' \ar[r] & X        }$$
in which the composition along the bottom row is the identity on $X$.
\end{proposition}

\begin{remark}
The idea is that $f$ is a retract of $i$ in a strong sense where
the sources of $f$ and $i$ are the same.
\end{remark}

\begin{proof}
Suppose that $i$ satisfies the conditions of the proposition.  Then $i$
is a cofibration by assumption, and it is a weak equivalence by Corollary
\ref{cor:n-equiv}.  Therefore, $f$ is a retract of a trivial cofibration,
so it is also a trivial cofibration.

Now suppose that $f$ is a trivial cofibration.
Use Proposition \ref{prop:factor-trcofib-fib} to factor $f$ as

$$\xymatrix{ A \ar[r]^i & X' \ar[r]^q & X,}$$
where $i$ satisfies the conditions of the proposition and $q$ is
a strong fibration.
By Proposition \ref{prop:lift-trcofib-strong-fib}, $f$ has the left
lifting property with respect to $q$.  The retract argument \cite[8.2.2]{PH}
finishes the proof.
\end{proof}

\begin{remark}
Weak equivalences of ind-spaces are slightly more difficult to define because
the language of local systems does not work so well.  See Remark
\ref{rem:ind-we} for a rough description of weak equivalences of ind-spaces.  
The model structure for ind-spaces is again a localization of the ``strict''
structure of Edwards and Hastings \cite[3.8]{EH}, so levelwise weak
equivalences are weak equivalences in our sense.

Fibrations of ind-spaces are levelwise fibrations.  Cofibrations are defined
by the appropriate left lifting property.  A dual notion of strong cofibration
is important in this context.
The technical heart of the theory lies in showing that
strong cofibrations have the left lifting property with respect to 
trivial fibrations.
\end{remark}

\section{Whitehead Theorem} 
\label{sctn:Whitehead}

Definition \ref{defn:we} is useful for studying the model structure,
but there are
other characterizations of weak equivalences.  We need a 
few preliminary definitions.

\begin{definition} 
\label{defn:cohlgy}
Let $X$ be a pro-space and let $L$ be a local system on $(X,\phi)$
represented by a functor $L_s: \Pi X_s \map \Ab$.
Define the {\em twisted cohomology groups} 

$$H^i(X; L) = \colim_{t \geq s} H^i (X_s; \phi^*_{ts} L_s).$$
If $f:Y \map X$ is a level cofibration, then

$$H^i(X,Y;L) = \colim_{t \geq s} H^i(X_t, Y_t; \phi^*_{ts} L_s).$$
\end{definition}

Note that there is a long exact sequence

$$\xymatrix@-2ex{
\cdots \ar[r] & H^i(X,Y; L) \ar[r] & H^i(X; L) \ar[r] & H^i(Y; f^*L)
   \ar[r] & H^{i+1}(X,Y; L) \ar[r] & \cdots    }$$
since filtered colimits are exact.

Recall that the Moore-Postnikov system \cite[8.9]{JPM} of a space $X$ is
a tower
$$ X \map \cdots \map P_n X \map P_{n-1} \map \cdots P_1 X \map P_0 X$$
such that $X \cong \lim_n P_n X$, each map $P_n X \map P_{n-1} X$ is 
a fibration, and each space $P_n X$ has the same homotopy groups as $X$ in
dimensions less than or equal to $n$ 
but has trivial homotopy groups in dimensions
greater than $n$.

\begin{definition} 
\label{defn:Postnikov}
If the functor $X: I \map \SSet$
is a pro-space, then $PX$ is the pro-space
given by the functor $\N \times I \map \SSet$ which sends $(n,s)$ to 
$P_n X_s$.
\end{definition}

Note that Artin and Mazur \cite{AM} used a similar construction $X^\natural$
instead of $PX$.

\begin{theorem} 
\label{thm:tfae}
Suppose that $f: X \map Y$ is a map of pro-spaces.  The following
conditions are equivalent:

(a) $f$ is a weak equivalence.

(b) $f$ is isomorphic to a level map of pro-spaces $g: Z \map W$ 
such that 
for all $n \geq 0$ and for all $s$, there exists some $t \geq s$ such that
for all basepoints $*$ in $Z_t$, there is a commutative diagram
$$\xymatrix{
\pi_n(Z_t,*) \ar[r]\ar[d] & \pi_n(W_t,*) \ar[d]\ar[dl] \\
\pi_n(Z_s,*) \ar[r]       & \pi_n(W_s,*).                    }$$

(c) $\pi_0(f)$ is an isomorphism of pro-sets, $\Pi_1 X \map f^*\Pi_1 Y$
is an isomorphism of pro-local systems on $X$, and for all $i$ and all
local systems $L$ on $Y$,
the map $H^i(Y; L) \map H^i(X; f^*L)$ is an isomorphism.

(d) $f$ is isomorphic to a level map $g: Z \map W$
indexed by a cofinite directed set $I$ for which there
is a strictly increasing function $n: I \map \N$ such that 
$g_s: Z_s \map W_s$ is an $n(s)$-equivalence.

(e) $P f$ is a strict weak equivalence.
\end{theorem}

We give the proof in Section \ref{sctn:cohlgy}.

Condition {\em (b)} is a less elegant but more technically convenient
definition of weak equivalences.
It is precisely Grossman's definition \cite{JG}
of weak equivalences for towers.
Condition {\em (c)} generalizes the well-known result about spaces that
a weak equivalence between simply connected spaces is the same
as a cohomology isomorphism.
Condition {\em (d)} gives the simplest description of weak equivalences.
This condition is usually easiest to use when working with the homotopy
theory of pro-spaces.
Condition {\em (e)} says that our homotopy theory is the 
``$P$-localization'' of the strict homotopy theory.
See Section \ref{sctn:strict} for the definition of strict weak
equivalences.

The proof that condition {\em (d)} is equivalent to
{\em (a)} relies
heavily on the model structure.  Condition {\em (d)} could have been the
definition of weak equivalences.  This would make the proofs of the model
structure axioms easier, but then it would be difficult to establish
a link between such
weak equivalences and pro-homotopy groups or cohomology groups.

\begin{remark} 
\label{rem:ind-we}
We take the dual to condition {\em (b)} as the definition of weak
equivalences of ind-spaces.  This definition is equivalent to a dual
version of condition {\em (d)}, but we do not know whether duals to
conditions {\em (a)}, {\em (c)}, and {\em (e)} exist.
\end{remark}

If $X$ is a pointed connected pro-space, then $\pi_n X$ is the pro-group
given by the diagram $\{ \pi_n X_s \}$.  These are the pro-homotopy
groups of Artin and Mazur \cite{AM}.  The following corollary indicates
that our weak equivalences appropriately generalize those of Artin and
Mazur.

\begin{corollary} 
\label{cor:pointed-we}
If $f: X \map Y$ is a map of connected pointed pro-spaces, then $f$ is a
weak equivalence if and only if $\pi_n f : \pi_n X \map \pi_n Y$ is an
isomorphism of pro-groups for all $n$.
\end{corollary}

\begin{proof}
We may assume that $f$ is a level map.
Suppose that the diagrams of condition {\em (b)}
of Theorem \ref{thm:tfae} exist.
Using Lemma \ref{lem:pro-iso}, it follows that 
$\pi_n f : \pi_n X \map \pi_n Y$ is an isomorphism of pro-groups.

Now suppose that $\pi_n f: \pi_n X \map \pi_n Y$ is an isomorphism of 
pro-groups.  Then the diagrams of Theorem \ref{thm:tfae} {\em (b)} exist 
for every $s$ when $*$ is the chosen basepoint of $X$.

Let $\sharp$ be any point of $X_t$.  Choose a path from $\sharp$
to $*$ in $X_t$.  This path induces an isomorphism between the diagram
of Theorem \ref{thm:tfae} {\em (b)} based at $*$ and the diagram based
at $\sharp$.
\end{proof}

\section{Comparison with Other Theories} 
\label{sctn:compare}

This section explains the difference between our homotopy category
of pro-spaces and the homotopy category of Artin and Mazur \cite{AM}.
We use the notations 
$$[ K, L ] = \Hom_{\Ho(\SSet)} ( K, L) $$
$$[X,Y]_{\pros} = \Hom_{\Ho(\pro\SSet)}(X,Y).$$
Let $c: \SSet \map \pro \SSet$ be the functor given by constant pro-spaces.

Artin and Mazur \cite{AM} defined a homotopy category of pro-spaces by
considering morphisms in $\lim_t \colim_s [X_s, Y_t]$.  Problems with
this category arise from the fact that diagrams commute only up
to homotopy.  More importantly, isomorphisms in their homotopy category
are {\em not} the same as maps that induce isomorphisms of pro-homotopy
groups.
  
If $X$ is cofibrant and $Y$ is fibrant, in general $[X,Y]_{\pros}$ is
not equal to $\lim_t \colim_s [X_s, Y_t]$.
In the pointed case,
there is a spectral sequence involving higher $\lim^i$ terms
relating the two sets.  In some special cases there is a simpler
relationship.

\begin{lemma} 
\label{lem:into-constant}
Let $X$ be a pro-space and let $K$ be a simplicial set with
finitely many nonzero homotopy groups.  Then 
$$[X,cK]_{\pros} \cong \colim_s [X_s, K].$$
\end{lemma}

\begin{proof}
Let $K'$ be a fibrant approximation of $K$.  Then $cK'$ is a fibrant
approximation of $cK$ in the category of pro-spaces.  Hence $[X,cK]_{\pros}$
equals $[X,cK']_{\pros}$ and $\colim_s [X_s, K]$ equals $\colim_s [X_s, K']$.
Therefore, it suffices to assume that $K$ is fibrant.

According to 
Definition \ref{defn:simp}, $\Map(X, cK) = \colim_s \Map(X_s, K)$.  Since
$\pi_0$ commutes with filtered colimits,
$$\pi_0 \Map(X, cK) \cong \colim_s \pi_0 \Map(X_s, K).$$
Since $X$ is cofibrant and $cK$ is fibrant, $[X, cK]_{\pros}$ is equal to
$\pi_0 \Map(X, cK)$, and $[X_s, K]$ is equal to $\pi_0 \Map(X_s, K)$ for
every $s$.  Therefore,
$$[X,cK]_{\pros} \cong \colim_s [X_s, K].$$
\end{proof}

For example, if $K = K(\pi,n)$, then

$$[X, K]_{\pros} \cong \colim_s [X_s,K] \cong \colim_s H^n(X; \pi)
   \cong H^n(X; \pi).$$
This explains why the ideas of Artin and Mazur were good enough for
studying ordinary cohomology.  However, other theories such as $K$-theory
present complications because they are not represented by spaces with
only finitely many nonzero homotopy groups.  Therefore, Artin and Mazur
were not able to define $K$-theory for pro-spaces.  We propose that $K$-theory
be defined as the functor represented by $BU$.  In order to calculate
$K$-theory, one must use the Postnikov system $PBU$ of $BU$.
This tower is a 
fibrant replacement for $BU$.  For a general pointed pro-space $X$,
there is a short exact sequence

$$0 \map \lim\nolimits^1_n \colim\nolimits_s [\Sigma X_s, P_n BU]
    \map [ X, BU]_{\pros}
    \map \lim\nolimits_n \colim\nolimits_s [X_s, P_n BU] \map 0.$$
When $X$ is the \'etale homotopy type of a scheme $S$ over a separably
closed field,
the third term in this sequence is Friedlander's definition \cite{etK-1}
of the \'etale $K$-theory of $S$.  Hence our definition differs by
a $\lim^1$ error term.

The theory of generalized cohomology for pro-spaces should clarify
the definition of the \'etale
$K$-theory of schemes.  Most likely the constructions of Dwyer and
Friedlander \cite{etK} can be put into
this framework for a cleaner description of the ideas.  We plan
to elaborate on this relationship in future work.

\begin{proposition} 
\label{prop:c-lim-adjoint}
The functors $c: \SSet \map \pro \SSet$ and
$\lim: \pro \SSet \map \SSet$ are a Quillen adjoint pair.
\end{proposition}

\begin{proof}
Let $X$ be a space and let $Y$ be a pro-space.  Then

$$\Hom_{\SSet}(X, \lim_s Y_s) = \lim_s \Hom_{\SSet}(X, Y_s) =
  \Hom_{\pro \SSet}(cX, Y).$$
Therefore $c$ and $\lim$ are an adjoint pair.

Since $c$ preserves cofibrations and trivial cofibrations, the functors
are a Quillen pair.
\end{proof}

\begin{corollary} 
\label{cor:c-lim-adjoint}
The functors $c$ and $\lim$ induce adjoint functors $\Ho(\SSet) \map
\Ho(\pro \SSet)$ and $\Ho(\pro \SSet) \map \Ho(\SSet)$.
\end{corollary}

The following proposition is a kind of dual to Lemma \ref{lem:into-constant}.

\begin{proposition} 
\label{prop:out-constant}
Let $K$ be a simplicial set, and let $Y$ be a pro-space.  Then
$$[cK, Y]_{\pros} = [K, \holim_s Y_s].$$
\end{proposition}

\begin{proof}
From Corollary \ref{cor:c-lim-adjoint}, we know that $[cK, Y]_{\pros}$
is equal to $[K, \lim_t Z_t]$, where $Z$ is a fibrant approximation of
$Y$.  It suffices to show that $\lim_t Z_t$ is weakly equivalent to 
$\holim_s Y_s$.

For each $s$, let $Y'_s$ be a functorial fibrant approximation of $Y_s$.
Then $Y \map Y'$ is a strict weak equivalence.  Now let $Z$ be a strictly
fibrant approximation of $PY'$.
Note that $Z$ is a fibrant approximation of $Y$.

Recall that $Y'_s$ is isomorphic to $\lim_n P_n Y'_s$ for each $s$.  Also,
$Y'_s$ is weakly equivalent to $\holim_n P_n Y'_s$ since each map
$P_n Y'_s \map P_{n-1} Y'_s$ is a fibration.

Now $\holim_s Y_s$ is weakly equivalent to $\holim_s Y'_s$ since $Y$ and $Y'$
are levelwise weakly equivalent.  Furthermore,
$\holim_s Y'_s$ is weakly equivalent
to $\holim_s \holim_n P_n Y'_s$.  In other words, $\holim_s Y_s$ is weakly
equivalent to $\holim_{n,s} PY'$.  Since $Z$ is a strictly fibrant 
approximation of $PY'$, $\lim_t Z_t$ is weakly equivalent to
$\holim_{n,s} PY'$ \cite[4.2]{EH}.
\end{proof}

Edwards and Hastings \cite[4.2]{EH} showed that homotopy limit is the derived
functor of limit with respect to the strict structure.  The proposition makes
the stronger claim that homotopy limit is the derived functor of limit
with respect also to our model structure.  In other words, pro-spaces with
isomorphic pro-homotopy groups have weakly equivalent homotopy limits,
even if they are not levelwise weakly equivalent.

The proposition implies that the pro-$K$-theory described above restricts
to ordinary topological $K$-theory on constant pro-spaces
because $BU$ is weakly equivalent to $\holim_n P_n BU$.

\section{Pro-Finite Completion} 
\label{sctn:Morel}

Now we make a comparison between the homotopy category of pro-spaces
and Morel's homologically defined homotopy theory on pro-simplicial
finite sets \cite{FM}.

Let $\SFS$ be the category of simplicial finite sets.  Thus, $\SFS$
consists of simplicial objects in the category of finite sets.  Note that
this is not the same as the category of finite simplicial sets since 
simplicial finite sets may have infinitely many non-degenerate simplices.

\begin{definition} \cite{FM} 
\label{defn:Morel}
Fix a prime number $p$.  A weak equivalence
in $\pro \SFS$ is a map $X \map Y$ such that $H^*(Y; \Z/p) \map H^*(X; \Z/p)$
is an isomorphism.  A cofibration in
$\pro \SFS$ is a map $X \map Y$ such that 
$\lim X \map \lim Y$ is a cofibration of spaces.  
Finally, a fibration in $\pro \SFS$ is a map with
the right lifting property with respect to trivial cofibrations.
\end{definition}

\begin{theorem} \cite{FM} 
\label{thm:Morel}
The above definitions give a model structure on the category $\pro \SFS$.
\end{theorem}

Pro-finite completion gives an adjoint pair of functors
between $\pro \SSet$ and $\pro \SFS$.
Sullivan \cite{DS} defined a topological pro-finite completion functor
$\wedge: \SSet \map \pro \SFS$.  It is 
the extension to simplicial objects of the functor from
sets to pro-finite sets sending $S$ to the pro-system of finite quotients
of $S$.  Pro-finite completion naturally extends to a functor
$\wedge: \pro \SSet \map \pro \SFS$.

The following lemma generalizes the usual adjointness property
of pro-finite completion.

\begin{lemma} 
\label{lem:pro-finite-adjoint}
Pro-finite completion $\wedge: \pro \SSet \map \pro \SFS$
is left adjoint to the
inclusion functor $i: \pro \SFS \map \pro \SSet$.
\end{lemma}

\begin{proof}
Let $X$ be a pro-space, and let $Y$ belong to $\pro \SFS$.  Then

\begin{eqnarray*}
\Hom_{\pro \SFS}(X^\wedge, Y) & = &
         \lim_s \colim_t \colim_u \Hom_{\SFS} ((X_t^\wedge)_u, Y_s) \\
 & = & \lim_s \colim_t \Hom_{\SSet} (X_t, Y_s) \\
 & = & \Hom_{\pro \SSet} (X, iY).      
\end{eqnarray*}
The second equality comes from the usual adjointness of pro-finite
completion, which says that $\colim_s\Hom_{\SSet} ((Z^\wedge)_s, W)
= \Hom_{\SSet} (Z, W)$
for any simplicial set $Z$ and any simplicial finite set $W$.
\end{proof}

\begin{proposition} 
\label{prop:Quillen-adjoint}
The functors $\wedge$ and $i$ of Lemma \ref{lem:pro-finite-adjoint}
form a Quillen pair of adjoint functors.
\end{proposition}

\begin{proof}
We must show that the left adjoint $\wedge$ preserves 
cofibrations and trivial cofibrations.

First consider an inclusion $f: X \map Y$
of simplicial sets.  The map 
$\lim_s (X^\wedge)_s \map \lim_t (Y^\wedge)_t$ is an inclusion.
This follows from the fact that if $x$ and $x'$ are two distinct 
$n$-simplices of $(X^\wedge)_s$, then there is a finite
quotient $(Y^\wedge)_t$ of $Y$ such that $f(x)$ and $f(x')$ are not equal.

Now consider an arbitrary cofibration of pro-spaces 
$f: X \map Y$.  Without loss of generality, we may assume that $f$ is
a levelwise cofibration.  
The map $\lim_t (X_s^\wedge)_t \map \lim_t (Y_s^\wedge)_t$
is an inclusion for each $s$, so the map

$$\lim_s \lim_t (X_s^\wedge)_t \map \lim_s \lim_t (Y_s^\wedge)_t$$
is also an inclusion.  
In other words, $X^\wedge \map Y^\wedge$ is a cofibration in $\pro \SFS$.

Next consider a trivial cofibration of pro-spaces
$f: X \map Y$.  We may assume that $f$ is
a levelwise cofibration.  Note that $H^*(Y; \Z/p) \map
H^*(X; \Z/p)$ is an isomorphism by Theorem \ref{thm:tfae}.

Morel \cite[1.2.2]{FM} observed
that $H^*(X_s; \Z/p) \cong H^*(X_s^\wedge; \Z/p)$ since $\Z/p$ is a finite
abelian group.  Therefore

\begin{eqnarray*}
H^*(X^\wedge; \Z/p) & = & \colim_s \colim_t H^*((X_s^\wedge)_t; \Z/p) \\
 & = & \colim_s H^*(X_s^\wedge; \Z/p) \\
 & \cong & \colim_s H^*(X_s; \Z/p) \\
 & = & H^*(X; \Z/p).
\end{eqnarray*}
Similarly, $H^*(Y^\wedge; \Z/p) \cong H^*(Y; \Z/p)$.  Hence
$H^*(Y^\wedge; \Z/p) \map H^*(X^\wedge; \Z/p)$ is isomorphic to the
map $H^*(Y; \Z/p) \map H^*(X; \Z/p)$, so it is also an isomorphism.
\end{proof}

\begin{corollary} 
\label{cor:Quillen-adjoint}
The functors $\wedge$ and $i$ of Lemma \ref{lem:pro-finite-adjoint}
induce an adjoint pair on the homotopy categories
$\Ho(\pro \SSet)$ and $\Ho(\pro \SFS)$.
\end{corollary}

\section{Strict Model Structure} 
\label{sctn:strict}

Edwards and Hastings \cite[3.3]{EH} defined a model structure on the category
of pro-spaces and used it to study shape theory.  We call this the
``strict'' structure because the weak equivalences are defined
levelwise.  We review their definitions and results.

\begin{definition} 
\label{defn:strict}
A map of pro-spaces is a {\em strict
weak equivalence} if it is isomorphic to
a levelwise weak equivalence.  A map of pro-spaces is a {\em strict
cofibration} if it is a cofibration.  A map of pro-spaces
is a {\em strict fibration} if it has the right lifting property with
respect to all strictly trivial cofibrations.
\end{definition}

\begin{definition} 
\label{defn:strict-strong-fib}
A map is a {\em strong strict fibration} if it is isomorphic to a
level map of pro-spaces $X \map Y$ indexed by a cofinite
directed set such that for all $t$,
$$X_t \map Y_t \times_{\lim\limits_{s<t} Y_s} \lim\limits_{s<t} X_s$$
is a fibration.
\end{definition}

\begin{theorem} \cite{EH}
The above definitions give a model structure on the category
$\pro \SSet$.
\end{theorem}

We call this the strict model structure on $\pro \SSet$.  Note that
Edwards and Hastings gave a more complicated definition for strict
weak equivalences.  The following proposition suffices to show that
their definition and our definition coincide.

\begin{proposition} 
\label{prop:compose-strict}
Strict weak equivalences are closed under composition.
\end{proposition}

\begin{proof}
Assume there are levelwise weak equivalences
$f: X \map Y$ and $g: Z \map W$ with an isomorphism $h: Y \map Z$.
We must construct a levelwise weak equivalence isomorphic to the 
composition $ghf$.

By adding isomorphisms to the systems for $f$ or $g$, make the
cardinalities of the index sets equal.
Choose an arbitrary isomorphism $\alpha$ from the index set of $g$ to
the index set of $f$.  

Define a function $t(s)$ inductively on height satisfying several
conditions.  Choose $t(s)$ large enough so that $Y_{t(s)} \map Z_s$
represents $h$, $t(s) \geq \alpha(s)$, and $t(u) < t(s)$ for all $u < s$.

The function $t$ defines cofinal subsystems of $X$ and $Y$.  Hence we
may assume that
$$\xymatrix{X \ar[r]^f & Y \ar[r]^h_\cong & Z \ar[r]^g & W}$$
is a level diagram indexed by a cofinite directed set $I$ 
where $f$ and $g$ are levelwise weak equivalences.
However, the composition is not necessarily a levelwise weak equivalence
because the map $h$ is not a levelwise weak equivalence.

Since $h$ is an isomorphism of pro-spaces, for every $s$, there exists 
$ t > s$ and a commutative diagram

$$\xymatrix{
Y_t \ar[r]\ar[d] & Z_t \ar[dl]\ar[d] \\
Y_s \ar[r]       & Z_s.                        }$$
By restricting to cofinal subsystems, we may assume that such a diagram
exists for every $t > s$.

Let $J$ be the directed set of indecomposable arrows of $I$.  The domain and
range functors $J \map I$ are both cofinal since $I$ is cofinite.
For each $\phi: t \map s$ in $J$, factor the map $Z_t \map Y_s$ as

$$\xymatrix{Z_t \ar[r]^{i_\phi} & A_\phi \ar[r]^{p_\phi} & Y_s,}$$
where $i_\phi$ is a cofibration and $p_\phi$ is a fibration.

Let $B_\phi$ be the pullback $X_s \times_{Y_s} A_\phi$, and let
$C_\phi$ be the pushout $W_t \amalg_{Z_t} A_\phi$.  These objects fit into
a commutative diagram

$$\xymatrix{
X_t \ar[r]\ar[d] & Y_t \ar[rr]\ar[dd]|(.5)\hole & &
        Z_t \ar[r]\ar[dd]|(.5)\hole \ar[dl] & W_t \ar[d]  \\
B_\phi \ar[rr]\ar[d] & & A_\phi \ar[rr]\ar[dl] & & C_\phi \ar[d] \\
X_s \ar[r] & Y_s \ar[rr] & & Z_s \ar[r]& W_s.   }$$

Note that $B_\phi \map A_\phi$ is a weak equivalence because it is a 
pullback of a weak equivalence along a fibration.  Also, 
$A_\phi \map C_\phi$ is a weak equivalence because it is a pushout of
a weak equivalence along a cofibration.  Hence the composition
$B \map A \map C$ is a levelwise weak equivalence.  This composition
is isomorphic to $ghf$ since $B \cong X$, $Y \cong A \cong Z$, and
$W \cong C$.
\end{proof}

The above proof works for any pro-category $\pro \C$ provided that $\C$ is
a proper model category.  In fact, a minor variation of the proof works
when $\C$ is either left proper or right proper.  
It is possible to prove the other parts of the
two-out-of-three axiom for strict weak equivalences
with similar techniques.  

\begin{proposition} 
\label{prop:strict-fib} \cite{EH}
A map of pro-spaces is a strict fibration if and only if it is a retract of
a strong strict fibration.
\end{proposition}

The relationship between the strict model structure and our model structure
is expressed in the following results.

Let $\Ho_{\strict}(\pro \SSet)$ be the homotopy category associated
to the strict structure.

\begin{proposition} 
\label{prop:localization}
The category $\Ho(\pro \SSet)$ is a localization of
$\Ho_{\strict}(\pro \SSet)$.
\end{proposition}

\begin{proof}
Every levelwise weak equivalence 
is a weak equivalence in the sense of Definition \ref{defn:we}.
\end{proof}

\begin{corollary} 
\label{cor:strict-fib}
If $f$ is a fibration, then $f$ is also a strict fibration.
\end{corollary}

\begin{proof}
The class of trivial cofibrations contains the class of strictly trivial
cofibrations by Proposition \ref{prop:localization}.
\end{proof}

\begin{corollary} 
\label{cor:strict-trfib}
A map of pro-spaces is a trivial fibration if and only if it is a strictly
trivial fibration.
\end{corollary}

\begin{proof}
Cofibrations are the same as strict cofibrations.
\end{proof}

\begin{proposition} 
\label{prop:compare-strict}
Let $X$ be a pro-space, and let $Y$ be a pro-space such that each $Y_s$ has
only finitely many non-zero homotopy groups.  Then

$$ [ X, Y]_{\pros} \cong \Hom_{\Ho_{\strict}(\pro \SSet)} (X, Y).$$
\end{proposition}

\begin{proof}
The condition on $Y$ ensures that its strictly fibrant replacement $Y'$
is also a fibrant replacement.  To calculate morphisms from $X$ to $Y$ in
either homotopy category, consider morphisms from $X$ to $Y'$ modulo
the simplicial homotopy relation.  Hence the morphisms are the same.
\end{proof}

For example, the proposition applies when $Y$ is a system of
Eilenberg-Maclane spaces.  
Edwards and Hastings described a relationship between 
homological algebra and the strict homotopy theory of such pro-spaces
\cite[Section 4]{EH}.
It follows that the relationship works just as well for our 
homotopy theory of pro-spaces.

\section{Limits and Colimits} 
\label{sctn:colim}

The rest of the paper concentrates on the technical details of the main
theorems stated in Section \ref{sctn:structure}.

We provide specific constructions of limits and colimits in
pro-categories.  The existence of all colimits seems to be a little-known
fact.
In this section only, consider pro-objects indexed by arbitrary
cofiltering categories, not just ones indexed by cofinite directed sets.

\begin{proposition} 
\label{prop:complete}
If $\C$ is complete, then $\pro \C$ is also complete.
If $\C$ is cocomplete, then $\pro \C$ is also cocomplete.
\end{proposition}

\begin{proof}
Artin and Mazur \cite[Appendix 4.2]{AM} showed that
$\pro \C$ has all equalizers (resp. coequalizers)
provided that $\C$ does.  It suffices to
construct arbitrary products (resp. coproducts) when these exist in $\C$.

Let $A$ be a set and let $X^\alpha$ be a pro-object for each $\alpha$ in $A$.
Let $I_\alpha$ be the cofiltering index category for $X^\alpha$.
Define $X = \prod_{\alpha \in A} X^\alpha$
to be the cofiltering system
with objects $\prod_{\alpha \in B} X^\alpha_{s_\alpha}$
and index category $I$
consisting of pairs $(B, (s_{\alpha}))$ where $B$ is a finite subset of $A$
and $(s_\alpha)$ is an element of
$\prod_{\alpha \in B} I_\alpha$.  A morphism
of $I$ from $(B, (s_\alpha))$ to $(C, (t_\alpha))$ corresponds to an
inclusion $C \subseteq B$ and morphisms
$s_\alpha \map t_\alpha$ in $I_\alpha$ for all $\alpha$ in $C$.
Use of finite subsets $B$ of $A$ is essential because we use
the fact that finite products commute with filtered colimits.

Direct calculation shows that for any $Y$ in $\pro \C$,
$$\Hom_{\pro \C}(Y, X) \cong
     \prod_{\alpha \in A} \Hom_{\pro \C}(Y, X^\alpha).$$
Thus arbitrary products exist.

To construct the coproduct, define $X = \coprod_{\alpha \in A} X^{\alpha}$
to be the cofiltering
system with objects $\coprod_{\alpha \in A} X^\alpha_{s_\alpha}$
and index
category $\{ (s_\alpha) \in \prod_{\alpha \in A} I_\alpha\}$.
Note that $\prod I_\alpha$ is cofiltering since each $I_\alpha$ is.

In order to show that $X$ has the correct universal mapping property, it
suffices to see that 

$$\Hom_{\pro \C}(X, Y) \cong
   \prod_{\alpha \in A} \Hom_{\pro \C}(X^\alpha, Y).$$
for $Y$ any object of
$\C$ ({\em i.e.}, $Y$ is a constant system in $\pro \C$).  This can be
checked directly, using the fact that colimits indexed on
product categories commute with the relevant products.
Thus arbitrary coproducts exist.
\end{proof}

\section{Retract Axioms} 
\label{sctn:retract}

The class of fibrations is obviously closed under retracts.  The class
of weak equivalences is closed under retracts because weak equivalences
are defined in terms of isomorphisms of pro-local systems.  We must show that 
the class of
cofibrations is also closed under retracts.  We prove a
general result and then apply it to cofibrations.

\begin{proposition} 
\label{prop:retract}
Let $\C$ be a category, and let $C$ be any class of maps in $\C$.  Define
the class $D$ as those maps in $\pro \C$ that are isomorphic to a level
map that belongs to $C$ levelwise.  Then $D$ is closed under retracts.
\end{proposition}

\begin{proof}
Suppose that $f: W \map Z$ is a retract of $g: X \map Y$, where
$g$ is a 
level map that belongs to $C$ levelwise.  Hence there is a commutative
diagram
$$\xymatrix{
W \ar[r] \ar[d]_f & X \ar[r] \ar[d]_g & W \ar[d]^f \\
Z \ar[r]          & Y \ar[r]          & Z          }$$
where the horizontal compositions are identity maps.
We must show that $f$ also
belongs to $D$.  Choose a level representative for $f$.

By adding isomorphisms to the systems for $f$ or $g$, make the
cardinalities of the index sets equal.  
Choose an arbitrary isomorphism $\alpha$ from the index set of $f$ to
the index set of $g$.

Define a function $t(s)$ inductively on height satisfying several
conditions.  First, choose $t(s)$ large enough so that $X_{t(s)} \map
W_s$ and $Y_{t(s)} \map Z_s$ represent respectively the maps $X \map W$
and $Y \map Z$.
Also,
choose $t(s)$ large enough so that $t(s) \geq \alpha(s)$.  Finally, choose
$t(s)$ large enough so that
for all $u < s$, $t(u) < t(s)$ with a commuting diagram

$$\xymatrix@-2ex{
 & X_{t(s)} \ar[rr] \ar[dl] \ar[dd]|!{[dl];[d]}\hole & & W_s \ar[dd]\ar[dl]\\
X_{t(u)} \ar[dd] \ar[rr] & & W_u \ar[dd] \\
 & Y_{t(s)} \ar[rr]|!{[ur];[dr]}\hole \ar[dl] & & Z_s \ar[dl] \\
Y_{t(u)} \ar[rr] & & Z_u.    }$$

Now the function $t$ defines cofinal subsystems
$\ol{X}$ and $\ol{Y}$ of $X$ and $Y$ where $\ol{X}$ and $\ol{Y}$ have the
same index sets as $W$ and $Z$.

Repeat this process on $\ol{g}: \ol{X} \map \ol{Y}$ to obtain another
function $u$ inducing cofinal subsystems $\ol{W}$ and $\ol{Z}$ of
$W$ and $Z$.  The result is a level diagram

$$\xymatrix{
\ol{W} \ar[r] \ar[d]_{\ol{f}} & \ol{X} \ar[d]_{\ol{g}} \ar[r] & W \ar[d]^f \\
\ol{Z} \ar[r] & \ol{Y} \ar[r] & Z          }$$
representing (up to isomorphism) $f$ as a retract of $g$.

Since $W \map X \map W$ and $Z \map Y \map Z$ are identity maps,
$u$ can be chosen so that the composites
$W_{u(s)} = \ol{W}_s \map \ol{X}_s \map W_s$
and $Z_{u(s)} = \ol{Z}_s \map \ol{Y}_s \map Z_s$ are structure maps of $W$
and $Z$ for all $s$.

Since $g$ belongs to $C$ levelwise,
the same is true for $\ol{g}$.  Define a pro-space $\hat{W}$
by starting with the system $W$ and replacing the single map
$W_{u(s)} \map W_s$ with the pair of maps $W_{u(s)} \map \ol{X}_s \map
W_s$. Define $\hat{Z}$ similarly.

Note that $W$ and $Z$ are cofinal subsystems of $\hat{W}$ and $\hat{Z}$
respectively,
so $W$ is isomorphic to $\hat{W}$ and $Z$ is isomorphic to $\hat{Z}$.
Thus it suffices to show that the level map 
$\hat{f}: \hat{W} \map \hat{Z}$ belongs to $D$.

The subsystem of $\hat{W}$ on objects $\{\ol{X}_s\}$ is 
also a cofinal subsystem, and the same is true for the subsystem on
objects $\{ \ol{Y}_s \}$ in $\hat{Z}$.  
Beware that the subsystem of $\hat{W}$
on $\{ \ol{X}_s \}$ is not isomorphic to $\ol{X}$ because the structure
maps are different.  The same warning applies to $\{ \ol{Y}_s \}$
and $\ol{Y}$.

Restrict
$\hat{f}$ to the subsystem $ \{ \ol{X}_s \} \map \{ \ol{Y}_s \}$.
This last map belongs to $C$ levelwise, so it belongs to $D$.
Hence $f$ also belongs to $D$ because $D$ is closed under isomorphisms.
\end{proof}

\begin{corollary} 
\label{cor:cofib-retract}
The class of cofibrations of pro-spaces is closed under retracts.
\end{corollary}

\begin{proof}
Apply Proposition \ref{prop:retract} to the class of all cofibrations
in $\SSet$.
\end{proof}

\section{Weak Equivalences} 
\label{sctn:we}

We begin with the two-out-of-three axiom.

\begin{proposition} 
\label{prop:2/3}
Let $f: X \map Y$ and $g: Y \map Z$ be maps of pro-spaces.
If any two of the maps $f$, $g$, and $gf$ are weak equivalences, then
so is the third.
\end{proposition}

\begin{proof}
For $n \geq 1$, the map
$\Pi_n X \map f^* g^* \Pi_n Z$ factors as 

$$\xymatrix@1{
\Pi_n X \map f^* \Pi_n Y \map f^* g^* \Pi_n Z.   }$$
Also, the map $\pi_0 X \map \pi_0 Z$ factors through $\pi_0 Y$.
This immediately proves two of the three cases.

For the third case, suppose that $f$ and $gf$ are weak equivalences.  Then 
$f^* \Pi_n Y \map f^* g^* \Pi_n Z$ is an isomorphism for all $n \geq 1$.
By Lemma \ref{lem:pro-pullback-pro-iso},
$\Pi_n Y \map g^* \Pi_n Z$ is also an isomorphism for all $n \geq 1$.
\end{proof}

The following lemma is a surprising 
generalization to pro-groups of an obvious fact about groups.  Bousfield
and Kan \cite[III.2.2]{BK} stated without proof a special case.

\begin{lemma} 
\label{lem:forget}
Let $U$ be the forgetful functor from pro-groups to pro-sets.
Then
a map $f$ of pro-groups is an isomorphism if and only if $U(f)$ is an
isomorphism of pro-sets.
\end{lemma}

\begin{proof}  
For simplicity write the group operations additively, even though the
groups are not necessarily abelian.

We may assume that $f:X \map Y$ is a level map.
If $f$ is an isomorphism, then $Uf$ is also an isomorphism by Lemma
\ref{lem:pro-iso}.

Now suppose that $Uf$ is an isomorphism.  By Lemma \ref{lem:pro-iso}
applied twice, for every $s$,
there exist $u \geq t \geq s$ and a commutative diagram

$$\xymatrix{
X_u \ar[r] \ar[d] & Y_u \ar[d]^\eta \ar[dl]_g \\
X_t \ar[r]^{f_t} \ar[d]_\phi & Y_t \ar[d] \ar[dl]_h \\
X_s \ar[r]        & Y_s   }$$
where the diagonal maps are not necessarily group homomorphisms.

However, the composite map $Y_u \map X_s$ is in 
fact a group homomorphism.
For every $x$ and $y$ in $Y_u$,
$$ h\eta(x + y) = h( \eta x + \eta y) = h( f_t g x + f_t g y) =
           hf_t( g x + g y)$$
because of commutativity in the top square.  Also,
$$ hf_t( g x + g y) = \phi( g x + g y) = \phi g x + \phi g y $$
because of commutativity in the bottom square.  Finally,
$\phi g = h \eta$.

Therefore, there is a commutative diagram of groups

$$\xymatrix{
X_u \ar[r] \ar[d] & Y_u \ar[d]\ar[dl] \\
X_s \ar[r]        & Y_s.              }$$
By Lemma \ref{lem:pro-iso}, $f$ is an isomorphism.
\end{proof}

The formal nature of Definition \ref{defn:we} is often too abstract for
comfort in technical situations.  The following proposition 
shows that conditions {\em (a)} and {\em (b)}
of Theorem \ref{thm:tfae} are equivalent, thus giving a less
natural but more concrete equivalent definition of weak equivalence.

\begin{proposition} 
\label{prop:we}
A level map of pro-spaces $f: (X,\phi) \map (Y,\eta)$ 
is a weak equivalence if and only if
for all $n \geq 0$ and for all $s$, there exists some $t \geq s$ such that
for all basepoints $*$ in $X_t$, there is a commutative diagram
$$\xymatrix{
\pi_n(X_t,*) \ar[r]\ar[d] & \pi_n(Y_t,*) \ar[d]\ar[dl] \\
\pi_n(X_s,*) \ar[r]       & \pi_n(Y_s,*).                    }$$
\end{proposition}

\begin{remark}
By the argument in the proof of 
Lemma \ref{lem:forget}, it is not important whether we assume that
the diagonal map is a group homomorphism or just a map of sets.
For convenience, we assume that it is a group homomorphism.
Note that the diagonal map is not geometrically induced; it just
exists abstractly.  The choice of $t$ may depend on $n$ and $s$, but it
must work for every basepoint.
\end{remark}

\begin{proof}
First suppose that $f$ is a weak equivalence.  Since $\pi_0 f$ is an
isomorphism, Lemma \ref{lem:pro-iso} gives the conclusion for $n = 0$.
For $n \geq 1$, Lemma \ref{lem:pro-iso} implies that,
for every $s$, there exists
a $t \geq s$ and a commutative diagram
$$\xymatrix{
\Pi_nX_t \ar[r]\ar[d] & f_t^* \Pi_n Y_t \ar[d] \\
\phi_{ts}^* \Pi_nX_s \ar[r]       & \phi_{ts}^* f_s^* \Pi_nY_s     }$$
of local systems on $X_t$.

In particular, for every basepoint $*$ in $X_t$, there is a commutative
diagram

$$\xymatrix{
\pi_n(X_t,*) \ar[r]\ar[d] & \pi_n(Y_t,*) \ar[d]\ar[dl] \\
\pi_n(X_s,*) \ar[r]       & \pi_n(Y_s,*).               }$$
This proves the ``only if'' part of the claim.

Now suppose that the diagrams in the statement of the proposition 
exist.  By Lemma \ref{lem:pro-iso},
$\pi_0 f$ is an isomorphism.  For $n \geq 1$, we use
the fact that a local system $L$ on a space $Z$ is determined up to
isomorphism by its value $L_x$ as a $\pi_1(Z,x)$-module
for one point $x$ in each component of $Z$.  For every $s$,
there exist $u \geq t \geq s$ such that for 
every basepoint $*$ in $X_u$ there are commutative diagrams

$$\xymatrix{
\pi_n(X_u,*) \ar[r]\ar[d] & \pi_n(Y_u,*) \ar[d]\ar[dl] &
   \pi_1(X_u,*) \ar[d]\ar[r] & \pi_1(Y_u,*) \ar[d]\ar[dl] \\
\pi_n(X_t,*) \ar[r]\ar[d] & \pi_n(Y_t,*) \ar[d]\ar[dl] &
   \pi_1(X_t,*) \ar[d]\ar[r] & \pi_1(Y_t,*) \ar[d]\ar[dl] \\
\pi_n(X_s,*) \ar[r] & \pi_n(Y_s,*) &
   \pi_1(X_s,*) \ar[r] & \pi_1(Y_s,*).                    }$$

A diagram chase like that in the proof of Lemma \ref{lem:forget}
shows that the map $\pi_n(Y_u,*) \map \pi_n(X_s,*)$ is
actually a map of $\pi_1(X_u,*)$-modules, even though the diagonal maps in
the left diagram above are not maps of $\pi_1(X_u,*)$-modules.  This
defines a commutative diagram

  $$\xymatrix{
\Pi_nX_u \ar[r]\ar[d] & f_u^*\Pi_nY_u \ar[d]\ar[dl] \\
\phi_{us}^*\Pi_nX_s \ar[r] & \phi_{us}^* f_s^* \Pi_nY_s        }$$
of local systems on $X_u$.

Hence $\Pi_nX \map f^*\Pi_nY$ is an isomorphism by Lemma \ref{lem:pro-iso}.
\end{proof}

\begin{corollary} 
\label{cor:n-equiv}
Suppose that $f: X \map Y$ is a level map of pro-spaces indexed by a
cofinite directed set
$I$ for which there is a strictly increasing function $n: I \map \N$ 
such that
$f_s: X_s \map Y_s$ is an $n(s)$-equivalence.  Then $f$ is a weak equivalence.
\end{corollary}

\begin{proof}
For any $s$ in $I$ and any $n \geq 0$, choose $t \geq s$ such that 
$n(t) > n$.  For every point $*$ in $X_t$, there is a commutative 
diagram of solid arrows 

$$\xymatrix{
\pi_n(X_t,*) \ar[r]^\cong \ar[d] & \pi_n(Y_t, *) \ar[d] \ar@{-->}[dl] \\
\pi_n(X_s,*) \ar[r] & \pi_n(Y_s, *).        }$$
Since the top horizontal map is an isomorphism, this diagram can
be extended to include the dashed arrow.
Thus $f$ satisfies the condition of Proposition \ref{prop:we}, so it
is a weak equivalence.
\end{proof}

\section{Fibrations} 
\label{sctn:fib}

The following lemma states some useful properties of strong fibrations that
follow directly from the definition.

\begin{lemma} 
\label{lem:fib}
If $f:X \map Y$ is a strong fibration, then for all $t$
the maps $f_t: X_t \map Y_t$ and
$g_t: \lim_{s<t}X_s \map \lim_{s<t}Y_s$
are co-$n(t)$-fibrations.
\end{lemma}

\begin{proof}
Let $h_t$ be the map 
$X_t \map Y_t \times_{\lim\limits_{s<t}Y_s} \lim\limits_{s<t}X_s$,
so $h_t$ is a co-$n(t)$-fibration by the definition of strong fibrations.

We use the lifting property characterization 
of co-$n(t)$-fibrations.  For the purposes of induction,
assume that the maps $f_s$ and
$g_s$ are co-$n(s)$-fibrations for all $s < t$.
In particular, they are all co-$n(t)$-fibrations since $n(s) \leq n(t)$.

Let $A \map B$ be an $n(t)$-cofibration.
A lift in the diagram

$$\xymatrix{
A \ar[r] \ar[d] & \lim_{s\leq u}X_s \ar[d]\ar[r] & \lim_{s<u}X_s \ar[d] \\
B \ar[r] \ar[urr]|!{[ur];[r]}\hole
     \ar@{-->}[ur] & \lim_{s\leq u}Y_s \ar[r] & \lim_{s<u}Y_s }$$
is the same as a lift for 

$$\xymatrix{
A \ar[r] \ar[d] & X_u \ar[d]^{h_u} \\
B \ar[r] \ar@{-->}[ur] & 
 Y_u \times_{\lim\limits_{s<u}Y_s} \lim\limits_{s<u}X_s.    }$$

A lift exists in the last diagram because $h_u$ is a co-$n(u)$-fibration.
Hence lifts
can be extended inductively to obtain lifts for $g_t$.
Thus $g_t$ is a co-$n(t)$-fibration.

The projection 
$p_t: Y_t \times_{\lim\limits_{s<t}Y_s} \lim\limits_{s<t}X_s \map Y_t$
is the pullback of $g_t$ along the base $Y_t \map \lim\limits_{s<t}Y_s$.
Therefore, it is also a co-$n(t)$-fibration since 
right lifting properties are preserved by pullbacks.
Now $f_t = p_t h_t$, so $f_t$ is a co-$n(t)$-fibration
since it is a composition of such maps.
\end{proof}

The follow proposition is useful for recognizing that particular maps
are strong fibrations.

\begin{proposition} 
\label{prop:fib}
A level map $f: X \map Y$ is a strong fibration if and only if for all
$t$, $f_t: X_t \map Y_t$ is a co-$n(t)$-equivalence and

$$h_t: X_t \map Y_t \times_{\lim\limits_{s<t}Y_s} \lim\limits_{s<t}X_s$$
is a fibration.
\end{proposition}

\begin{proof}
One direction was proved in Lemma \ref{lem:fib}.

Use the notation of Lemma \ref{lem:fib}.
Suppose that $f_t$ is a co-$n(t)$-equivalence and $h_t$ is a 
fibration for all $t$.  It suffices to show that $h_t$ is also a
co-$n(t)$-equivalence.  For the purposes of induction, assume
that $h_s$ is a co-$n(s)$-equivalence for $s < t$.

As in the proof of Lemma \ref{lem:fib}, $g_t$ is a co-$n(t)$-fibration
by the inductive assumption.  Hence $p_t$ also is a co-$n(t)$-fibration.
The maps $p_t$ and $f_t$ are  co-$n(t)$-equivalences.
Since $f_t = p_t h_t$, it follows that $h_t$ is also a co-$n(t)$-equivalence.
\end{proof}

Our next goal is to show that strong fibrations are also fibrations.
This does not follow immediately from the definitions.
First we need two preliminary lemmas about simplicial sets.

\begin{lemma} 
\label{lem:technical-lift}
Suppose given a commutative diagram of simplicial sets

$$\xymatrix{
A \ar[d]_j \ar[r] & E_0 \ar[d]_{p_0} \ar[r] & E_1 \ar[d]_{p_1} \ar[r] &
\cdots \ar[r] & E_N \ar[d]_{p_N} \ar[r] & E_{N+1} \ar[d]_{p_{N+1}} \\
X \ar[r]          & B_0 \ar[r]              & B_1 \ar[r]              & 
\cdots \ar[r]     & B_N \ar[r]              & B_{N+1}                  }$$
such that $j$ is a cofibration, every $p_i$ is a fibration, and 
$p_{N+1}$ is a co-$N$-fibration.

For every choice of basepoint $*$ in $E_i$, there is an induced map
$F_i \map F_{i+1}$ where $F_i$ and $F_{i+1}$ are fibers of $p_i$ and
$p_{i+1}$.  Suppose that for every $i \leq N-1$
and every point $*$ in $E_i$, the map
$\pi_i(F_i,*) \map \pi_i(F_{i+1},*)$ is trivial.
Finally, suppose that the image
of $\pi_0 E_0$ in $\pi_0 B_0$ contains the image of $\pi_0 X$.
Then there is a lift in the diagram

$$\xymatrix{
A \ar[r] \ar[d]  & E_{N+1} \ar[d] \\
X \ar[r] \ar@{-->}[ru] & B_{N+1}.       }$$
\end{lemma}

\begin{proof}
There exists a lift in the diagram

$$\xymatrix{
A \ar[r] \ar[d]  & E_0 \ar[d] \\
X^{(0)}  \ar[r] \ar@{-->}[ru] & B_0 }$$
because the image of $\pi_0 E_0$ in $\pi_0 B_0$ contains the image
of $\pi_0 X$.

Now consider a diagram

$$\xymatrix{
X^{(m)} \ar[r] \ar[d]  & E_m \ar[r] & E_{m+1} \ar[d] \\
X^{(m+1)} \ar[r]         & B_m \ar[r] & B_{m+1}.       }$$
We have already shown that this diagram exists for $m = 0$.  We show by
induction that this diagram exists for all $m \leq N$.  The induction
step amounts to finding a lift in the above diagram.

We construct a lift one simplex at a time.  Therefore, it suffices to
assume that $X^{(m)} = \partial \Delta^{m+1}$
and $X^{(m+1)} = \Delta^{m+1}$.  The obstruction
to finding a lift is an element of $\pi_m F_{m+1}$, and this obstruction
lies in the image of $\pi_m F_m$.  When $m < N$ the map 
$\pi_m F_m \map \pi_m F_{m+1}$ is trivial, so
there is no obstruction to lifting.

By induction on $m$, there is a diagram
$$\xymatrix{
X^{(N)} \ar[r]\ar[d] & E_N       \ar[r] & E_{N+1} \ar[d]  \\
X         \ar[r]       & B_N  \ar[r]    & B_{N+1}.      }$$
A lift exists in this diagram since
the left vertical arrow is an $N$-cofibration and
the right vertical arrow is a co-$N$-fibration.
The map $X \map E_{N+1}$ is the desired lifting.
\end{proof}

\begin{lemma} 
\label{lem:zero}
Suppose that $n \geq 0$ and given a commutative diagram of simplicial sets

$$\xymatrix{
E_2 \ar[r] \ar[d]_{p_2} & E_1 \ar[r] \ar[d]_{p_1} & E_0 \ar[d]_{p_0} \\
B_2 \ar[r]              & B_1 \ar[r]              & B_0            }$$
such that each $p_i$ is a fibration.  For every choice of
basepoint $*$ in $E_2$, there are induced maps
$F_2 \map F_1 \map F_0$ where $F_i$ is the fiber of $p_i$ 
over the basepoint.

Suppose that for every $*$ in $E_2$, the maps
$\pi_n(E_2,*) \map \pi_n(E_1,*)$ and $\pi_{n+1}(B_1,*) \map \pi_{n+1}(B_0,*)$
factor through $\pi_n(B_2,*)$ and $\pi_{n+1}(E_0,*)$ respectively.  Then
for every $*$, the map $\pi_n(F_2,*) \map \pi_n(F_0,*)$ is trivial.
\end{lemma}

\begin{proof}
Chase the commutative diagram with exact rows

$$\xymatrix{
&&\pi_n(F_2,*)\ar[r]\ar[d] & \pi_n(E_2,*)\ar[r]\ar[d] & \pi_n(B_2,*)\ar[dl] \\
& \pi_{n+1}(B_1,*) \ar[r]\ar[ld]\ar[d] & \pi_n(F_1,*) \ar[r]\ar[d] &
        \pi_n(E_1,*) \\
\pi_{n+1}(E_0,*) \ar[r] & \pi_{n+1}(B_0,*) \ar[r] & \pi_n(F_0,*).   }$$
\end{proof}

The next proposition shows that strong fibrations are also fibrations.  This 
result is the technical heart of the whole theory.

\begin{proposition} 
\label{prop:lift-trcofib-strong-fib}
Strong fibrations have the right lifting property with respect to
trivial cofibrations.
\end{proposition}

\begin{proof}
Consider a commutative diagram of pro-spaces

$$\xymatrix{
A \ar[r] \ar[d]_j & E \ar[d]^p \\
X \ar[r]        & B            }$$
such that $j$ is a trivial cofibration and $p$ is a strong fibration.
We may assume that $j$ is a level map that is a 
levelwise cofibration.  By
choosing an appropriate cofinal subsystem for $j$ as in the proof of
Proposition \ref{prop:retract}, we may additionally
assume that the square diagram is a level diagram.  This choice preserves
$j$ as a levelwise cofibration.

We construct a lifting by induction.  Suppose that for all $s < t$, there
exists $a(s) \geq s$ and a map $X_{a(s)} \map E_s$ such that
for all $u < s < t$,
$a(u) \leq a(s)$ and there is a commutative diagram

$$\xymatrix@-2ex{
 & E_s \ar[rr] \ar[dd]|!{[dl];[d]}\hole & & E_u \ar[dd] \\
A_{a(s)} \ar[rr] \ar[ur] \ar[dd] & & A_{a(u)} \ar[ur] \ar[dd] & \\
 & B_s \ar[rr]|(.5)\hole |(.667)\hole & & B_u \\
X_{a(s)} \ar[rr] \ar[ur] \ar[ruuu]|!{[uu];[uurr]}\hole &
                                         & X_{a(u)}. \ar[ur] \ar[ruuu] & }$$

To extend the map $X \map E$ to level $t$, we must 
find a lift in the diagram
$$\xymatrix{
A_{a(t)} \ar[r] \ar[d] & E_t \ar[d] \\
X_{a(t)} \ar[r] \ar@{-->}[ur] &
      B_t \times_{\lim\limits_{s<t}B_s} \lim\limits_{s<t}E_s  }$$
for some $a(t) \geq t$
such that $a(s) < a(t)$ for all $s < t$.

Write the map 
$E_t \map B_t \times_{\lim\limits_{s<t} B_s} \lim\limits_{s<t} E_s$
as $q: W \map Z$.  Because $p: E \map B$ is a strong fibration,
$q$ is a co-$N$-fibration for some $N$.

Choose $t_{2N}$ so that $t_{2N} > a(s)$ for all $s<t$.  Now select
$t_{2N-1}, t_{2N-2}, \ldots, t_1, t_0$ so that $t_i > t_{i+1}$
and there exist commutative diagrams

$$\xymatrix{
\pi_n(A_{t_i},*) \ar[r] \ar[d] & \pi_n(A_{t_{i+1}},*) \ar[d] \\
\pi_n(X_{t_i},*) \ar[r] \ar[ru] & \pi_n(X_{t_{i+1}},*). }$$
for all $0 \leq n \leq N$ and all basepoints $*$ in $A_{t_i}$.
This is possible since $j$ is a weak equivalence and there are only
finitely many conditions on the choice of each $t_i$.

Finally, choose $a(t)$ so that $a(t) > t_0$ and there exists a commutative
diagram

$$\xymatrix{
\pi_0 A_{a(t)} \ar[r] \ar[d] & \pi_0 A_{t_0} \ar[d] \\
\pi_0 X_{a(t)} \ar[r] \ar[ru] & \pi_0 X_{t_0}. }$$

Functorially factor each map $A_{t_i} \map X_{t_i}$ as
$\xymatrix@1{A_{t_i} \ar[r]^{a_i} & Y_{t_i} \ar[r]^{b_i} & X_{t_i}}$, where
$a_i$ is a trivial cofibration and $b_i$ is a fibration.

Choose a basepoint $*$ in $Y_{t_i}$.  Since $a_i$ is a weak equivalence,
there exists a basepoint $\sharp$ in $A_{t_i}$ such that there is a path in
$Y_{t_i}$ from $*$ to the image of $\sharp$.  For $0 \leq n \leq N$,
there is a diagram of solid arrows

$$\xymatrix@-2ex@C-3ex{
 & \pi_n(A_{t_i},\sharp) \ar[dl] \ar[rr]
      \ar[dd]|(.333)\hole|(.5)\hole & &
      \pi_n(A_{t_{i+1}},\sharp) \ar[dl] \ar[dd]                         \\
\pi_n(X_{t_i},\sharp) \ar[rr] \ar[dd] \ar[rrru] & &
      \pi_n(X_{t_{i+1}},\sharp) \ar[dd] \\
 & \pi_n(Y_{t_i},*) \ar[dl] \ar[rr]|!{[ur];[dr]}\hole & &
      \pi_n(Y_{t_{i+1}},*) \ar[dl]                                  \\
\pi_n(X_{t_i},*) \ar[rr] \ar@{-->}[urrr]|!{[uurr];[rr]}\hole & &
    \pi_n(X_{t_{i+1}},*)     }$$
where the vertical arrows are induced by the choice of path.  Since the 
vertical arrows are all isomorphisms, the dashed arrow also exists.

For similar reasons, there is also a commutative diagram

$$\xymatrix{
\pi_0 A_{a(t)} \ar[r] \ar[d] & \pi_0 Y_{t_0} \ar[d] \\
\pi_0 X_{a(t)} \ar[r] \ar[ru] & \pi_0 X_{t_0}. }$$

Now we have

$$\xymatrix{
A_{a(t)}\ar[r]\ar[dd] & A_{t_0} \ar[r]\ar[d]_{a_0} &
     A_{t_1} \ar[r]\ar[d]_{a_1} & \cdots \ar[r] &
     A_{t_{2N}} \ar[d]_{a_{2N}}\ar[r] & W \ar[dd]^q \\
 & Y_{t_0} \ar[r]\ar[d]_{b_0} & Y_{t_1} \ar[r]\ar[d]_{b_1} & \cdots \ar[r] &
     Y_{t_{2N}} \ar[d]_{b_{2N}}\ar@{-->}[ur] \\
X_{a(t)}\ar[r] & X_{t_0} \ar[r] & X_{t_1} \ar[r] & \cdots \ar[r] &
     X_{t_{2N}} \ar[r] & Z }$$
where the dashed arrow exists because $a_{2N}$ is a trivial cofibration and
$q$ is a fibration.  So we need only find a lift for the diagram

$$\xymatrix{
A_{a(t)} \ar[r] \ar[d] & Y_{t_0} \ar[r] & Y_{t_2} \ar[r] & \cdots \ar[r] &
      Y_{t_{2N-2}} \ar[r] & Y_{t_{2N}} \ar[r] & W \ar[d] \\
X_{a(t)} \ar[r] & X_{t_0} \ar[r] & X_{t_2} \ar[r] & \cdots \ar[r] &
      X_{t_{2N-2}} \ar[r] & X_{t_{2N}} \ar[r] & Z.  }$$
Lemma \ref{lem:zero} tells us that this diagram satisfies the hypotheses of
Lemma \ref{lem:technical-lift}.  Hence the desired lift exists.
\end{proof}

\section{Lifting and Factorization Axioms} 
\label{sctn:lift-factor}

We now prove the lifting and factorizations axioms.

\begin{proposition} 
\label{prop:factor-trcofib-fib}
If $f:X \map Y$ is a map of pro-spaces, then $f$
factors (not functorially) as

$$\xymatrix{X \ar[r]^i & Z \ar[r]^p & Y,}$$
where $i$ is a trivial cofibration and $p$ is a strong fibration.
\end{proposition}

\begin{proof}
We may assume that $f$ is a level map.  We construct the factorization
inductively.

Assume for sake of induction that the factorization is already constructed on 
all indices less than $t$.
Recall the height function $h(t)$ from Definition \ref{defn:height}.
Factor
$$X_t \map Y_t \times_{\lim\limits_{s<t} Y_s} \lim\limits_{s<t} Z_s$$
as 
$$\xymatrix{X_t \ar[r]^{i_t} & Z_t \ar[r]^-{q_t} &
     Y_t \times_{\lim\limits_{s<t} Y_s} \lim\limits_{s<t} Z_s}$$
where $i_t$ is an $h(t)$-cofibration and $q_t$ is a co-$h(t)$-fibration.
Let $p_t: Z_t \map Y_t$ be
the projection map induced by $q_t$.  This extends the factorization
to level $t$.

Now $p$ satisfies the definition of strong fibration by construction.
Also, $i$ is a levelwise cofibration by construction.  Finally,
$i$ is a weak equivalence by Corollary \ref{cor:n-equiv}.
\end{proof}

\begin{proposition} 
\label{prop:factor-cofib-trfib}
If $f: X \map Y$ is a map of pro-spaces, then $f$
factors (not functorially) as 

$$\xymatrix{X \ar[r]^i & Z \ar[r]^p & Y,}$$
where $i$ is a cofibration and $p$ is a trivial fibration.
\end{proposition}

\begin{proof}
Use the strict model structure to
factor $f$ as a cofibration $i$ followed by a strictly trivial fibration $p$.
Then $p$ is also a trivial fibration by Corollary \ref{cor:strict-trfib}.
\end{proof}

\begin{remark}
We rely here on factorizations in the strict model structure.
These factorizations are constructed similarly to the method of the proof
of Proposition \ref{prop:factor-trcofib-fib}.
\end{remark}

One of the lifting axioms is automatic from the definition of fibrations.
The other axiom is proved in the following proposition.

\begin{proposition} 
\label{prop:lift-cofib-trfib}
Cofibrations have the left lifting property 
with respect to trivial fibrations.
\end{proposition}

\begin{proof}
Consider a commutative diagram of pro-spaces

$$\xymatrix{
A \ar[r] \ar[d]_j & E \ar[d]^p \\
X \ar[r]        & B            }$$
such that $j$ is a cofibration and $p$ is a trivial fibration.
We may assume that $j$ is a level map that is a 
levelwise cofibration.  By
choosing an appropriate cofinal subsystem for $j$, we may additionally
assume that the square diagram is a level diagram.  This choice preserves
$j$ as a levelwise cofibration.

Use the strict model structure to factor $p$ as 
$\xymatrix@1{E \ar[r]^i & Y \ar[r]^q & B}$,
where $i$ is a levelwise cofibration
and $q$ is a strictly trivial fibration.  Note that
$i$ is a weak equivalence by the two-out-of-three axiom as proved 
in Proposition \ref{prop:2/3}.

In the diagram

$$\xymatrix{
                           & E \ar[d]^i \ar[r]^=             & E \ar[d]^p \\
A \ar[d]_j \ar[ru] \ar[r]  & Y \ar[d]^q \ar[r] \ar@{-->}[ur] & B \\
X \ar[r] \ar@{-->}[ru]     & B, \ar[ru]_=                          }$$
the lift in the lower left square exists because of the strict model
structure, and the
lift in the upper right square exists by the definition of fibrations.
The composition $X \map Y \map E$ is the desired lift.
\end{proof}

\section{Simplicial Model Structure} 
\label{sctn:simp}

Recall that a model structure on a category $\C$ 
is simplicial if for every $X$ in $\C$ and every simplicial set $K$,
there are functorial constructions $X \otimes K$ (``tensor'')
and $X^K$ (``cotensor'') in $\C$ satisfying
certain associativity and unit conditions.  Also, for every $X$ and $Y$ in
$\C$, there is a simplicial function complex $\Map(X,Y)$.
These three constructions are related by the adjunctions

$$\Map(X \otimes K, Y) \cong \Map(K, \Map(X,Y)) \cong \Map(X, Y^K).$$
Finally, $\Map(-,-)$ must interact appropriately with the model structure as
follows.  If $i:A \map X$ is a cofibration in $\C$
and $p:E \map B$ is a fibration in $\C$, then

$$\Map(X,E) \map \Map(A,E) \times_{\Map(A,B)} \Map(X,B)$$
is a fibration that is a weak equivalence if either $i$ or $p$ is a 
weak equivalence.

We begin with a general proposition showing that tensors and cotensors
defined for {\em finite} simplicial sets automatically extend to all
simplicial sets.

\begin{proposition} 
\label{prop:extend}
Suppose that $\C$ is a model category with a simplicial function
complex $\Map(-, -)$.  Also suppose that the tensor $X \otimes K$
and the cotensor $X^K$ are
defined for all objects $X$ of $\C$ and all {\em finite} simplicial sets
$K$ so that the axioms for a simplicial model structure are
satisfied when they make sense.  Then
the definitions of tensor and cotensor can be
extended to provide a simplicial model structure for $\C$.
\end{proposition}

\begin{proof}
For any simplicial set $K$,
let $K^{\fin}$ be the filtering system of finite
subspaces of $K$.  For an object $X$ of $\C$, define $X \otimes K$ to
be $\colim_s (X \otimes K_s^{\fin})$ and $X^K$ to be $\lim_s X^{K_s^{\fin}}$.
Using the fact that the system $K^{\fin} \times L^{\fin}$ is cofinal
in the system $(K \times L)^{\fin}$, the required isomorphisms

$$X \otimes (K \times L) \cong (X \otimes K) \otimes L$$
and

$$\Map(X \otimes K, Y) \cong \Map(K, \Map(X, Y)) \cong \Map(X, Y^K)$$
can be verified directly.
\end{proof}

\begin{definition} 
\label{defn:simp}
If $X$ and $Y$ are pro-spaces and $K$ is a simplicial set, define
$$\Map(X,Y) = \Hom_{\pro \SSet}(X \times \Delta^{\bullet}, Y)
  = \lim_s \colim_t \Map(X_t, Y_s),$$
$$X \otimes K = \colim_s (X \times K_s^{\fin}),$$
and
$$Y^K = \lim_s (Y^{K_s^{\fin}}).$$
\end{definition}

For an arbitrary pro-space $X$ and a simplicial set
$K$, $X \times K$ can be constructed as the levelwise product with $K$.
Also, $\lim_s Y^{K_s^{\fin}}$ can be constructed as the system
$\{ Y_t^{K_s^{\fin}} \}$, indexed by all pairs $(s, t)$
in the product of the index categories.

Note that $X \otimes K$ is not in general isomorphic to $X \times K$
because finite limits do not always commute with filtered colimits
in the category of pro-spaces.  However, if $K$ is
finite, then $X \otimes K$ is isomorphic to $X \times K$ since $K$
itself is the terminal object of $K^{\fin}$.  Also, when $K$ is finite,
$Y^K$ is the system $\{Y_s^K\}$ with the same index
category as that of $Y$.

\begin{proposition} 
\label{prop:simp}
The above definitions make $\pro \SSet$ into a simplicial model category.
\end{proposition}

\begin{proof}
By Proposition \ref{prop:extend},
it suffices to check the axioms only for finite
simplicial sets.  
Most of the axioms are obvious; we verify only the
non-trivial ones here.

Let $X$ and $Y$ be arbitrary pro-spaces, and let $K$ be a finite
simplicial set.  We use the fact that $\Hom_{\SSet}(K, \colim_s Z_s)$
is equal to $\colim_s \Hom_{\SSet} (K, Z_s)$ for any 
filtered system $Z$
of simplicial sets because $K$ is finite.  It follows by direct
calculation that

$$\Map(X \otimes K, Y) \cong \Map (X, Y^K) \cong \Map (K, \Map(X,Y)).$$

We now show that the map
$$f: \Map(B,X) \map \Map(A,X) \times_{\Map(A,Y)} \Map(B,Y)$$
is a fibration whenever $i: A \map B$ is a cofibration and
$p: X \map Y$ is a fibration and that this map 
is a trivial fibration if either $i$ or $p$ is trivial.  We proceed
by showing that $f$ has the relevant right lifting property.  Let
$j: K \map L$ be a generating cofibration or a generating trivial
cofibration.  Note that $K$ and $L$ are finite simplicial sets.

By adjointness, it suffices to show that the map 

$$g: A \otimes L \amalg_{A \otimes K} B \otimes K \map B \otimes L$$
is a cofibration that is trivial if either $i$ or $j$ is trivial.
 
We may assume that $i$ is a levelwise cofibration.  For every 
$s$, $A_s \map B_s$ is a cofibration.  Therefore,
the map 
$$A_s \otimes L \amalg_{A_s \otimes K} B_s \otimes K \map B_s \otimes L$$
is also a cofibration.  This is a standard fact about 
simplicial sets.  Thus $g$ is a levelwise cofibration.

In order to show that $g$ is trivial whenever $i$ or $j$ is, it suffices
to show that the map $A \otimes K \map B \otimes K$ is trivial if $i$ is
trivial and that the map $A \otimes K \map A \otimes L$ is trivial
if $j$ is trivial.
This reduction follows from the two-out-of-three axiom,
the fact that trivial cofibrations are preserved by pushouts, and the
commutative diagram

$$\xymatrix{
A \otimes K \ar[r]\ar[d] & A \otimes L \ar[d]\ar[ddrr] \\
B \otimes K \ar[r]\ar[rrrd] &
 A \otimes L \amalg_{A \otimes K} B \otimes K \ar[drr]     \\
 & & & B \otimes L.    }$$

First suppose that $j$ is trivial.
The map $A \otimes K \map A \otimes L$ is a levelwise
weak equivalence, so it is a weak equivalence of pro-spaces.

Now suppose that $i$ is trivial.  
Since $A \otimes K \map B \otimes K$ is constructed by levelwise
product with $K$, condition {\em (b)} of Theorem \ref{thm:tfae} is
easily verified.
\end{proof}

\section{Properness} 
\label{sctn:proper}

We now show that the model structure of Theorem \ref{thm} is proper.
Recall that a model structure is left proper if weak equivalences
are preserved under pushout along cofibrations.  Dually, a model structure
is right proper if weak equivalences
are preserved under pullback along fibrations.

\begin{proposition} 
\label{prop:proper}
The simplicial model structure of Theorem \ref{thm} is left and right
proper.
\end{proposition}

\begin{proof}
Left properness follows immediately from the fact that all pro-spaces
are cofibrant.  We must show that the model structure is
right proper.

Let $p: E \map B$ be a fibration and let
$f: X \map B$ be a weak equivalence.  Use Theorem \ref{thm:tfae} to
suppose that
$p$ and $f$ are level maps with
the same cofinite directed index set $I$ for which 
there is a strictly
increasing function $n: I \map \N$ such that $f_s$ is a $n(s)$-equivalence.
Let $P$ be the pullback $X \times_B E$,
which is constructed levelwise.  We must show that
the projection $P \map E$ is a weak equivalence.

We start with a special case.  First suppose that $p$ is a levelwise
fibration.
Let $*$ be a basepoint in $P_s$.  This yields a diagram

$$\xymatrix{
F \ar[r] \ar[d]_= & P_s \ar[r] \ar[d] & X_s \ar[d]^{f_s} \\
F \ar[r]          & E_s \ar[r]_{p_s}  & B_s        }$$
in which the rows are fiber sequences.  From the $5$-lemma applied
to the long exact sequences of homotopy groups of the fibrations,
$P_s \map E_s$ is also an $n(s)$-equivalence.  By Theorem \ref{thm:tfae},
$P \map E$ is a weak equivalence.

Now let $p$ be an arbitrary fibration.
By Proposition \ref{prop:fib-retract}, there exists a strong fibration
$q:E' \map B$ such that $p$ is a retract of $q$.
Note that $q$ is a levelwise fibration by Lemma \ref{lem:fib}.

Consider the commutative diagram

$$\xymatrix@-2ex{
 & P \ar[dl] \ar[rr]\ar[dd]|(.5)\hole & & P' \ar[dl]\ar[rr]\ar[dd]|(.5)\hole
       & & P \ar[dl]\ar[dd] \\
E \ar[dd]\ar[rr] & & E' \ar[dd]\ar[rr] & & E \ar[dd] \\
 & X \ar[dl]\ar[rr]|(.5)\hole & & X \ar[dl]\ar[rr]|(.5)\hole & & X \ar[dl] \\
B \ar[rr] & & B \ar[rr] & & B,      }$$
where $P' = X \times_B E'$.  This diagram is a retract of squares in the
sense that all of the horizontal compositions are identity maps.  
The map $P' \map E'$ is a weak equivalence by the special case.  Since 
weak equivalences are closed under retracts, the map $P \map E$ is also
a weak equivalence.
\end{proof}

\section{Alternative Characterizations of Weak Equivalences}
\label{sctn:cohlgy}

We finish here the proof of Theorem \ref{thm:tfae} describing weak
equivalences in other terms.  For expository clarity, we split
the theorem into several parts.  The equivalence of {\em (a)} and {\em (b)}
was shown in Proposition \ref{prop:we}.

\begin{proposition} 
\label{prop:tfae1}
A map of pro-spaces is a weak equivalence if and only if
it is isomorphic to a level map $g: Z \map W$
indexed by a cofinite directed set $I$ for which there
is a strictly increasing function $n: I \map \N$ such that 
$g_s: Z_s \map W_s$ is an $n(s)$-equivalence.
\end{proposition}

\begin{proof}
Corollary \ref{cor:n-equiv} showed that a map $g$ satisfying the
conditions of the proposition is a weak equivalence.

Now suppose that $f$ is a weak equivalence.  We may assume that $f$ is a level
map.  Use Proposition \ref{prop:factor-trcofib-fib} to factor $f$ as 

$$\xymatrix{
X \ar[r]^i & Z \ar[r]^p & Y,  }$$
where $i$ is a trivial cofibration and $p$ is a trivial fibration.
By Corollary \ref{cor:strict-trfib}, $p$ is also a strictly trivial
fibration.  In particular, $p$ is isomorphic to a 
levelwise weak equivalence.  The proof
of Proposition \ref{prop:factor-trcofib-fib} indicates that $i$
satisfies the conditions of the proposition.  By an argument similar
to the proof of Proposition \ref{prop:compose-strict}, $f$ also satisfies
the conditions of the proposition.
\end{proof}

Recall the Moore-Postnikov functor $P$ from Definition \ref{defn:Postnikov}.

\begin{lemma} 
\label{lem:postnikov-we}
The canonical map $X \map P X$
is a weak equivalence for any pro-space $X$.
\end{lemma}

\begin{proof}
Condition {\em (b)} of Theorem \ref{thm:tfae} is easily verified.
\end{proof}

\begin{proposition} 
\label{prop:tfae2}
A map of pro-spaces $f: X \map Y$ is a weak equivalence if and only if
$P f$ is a strict weak equivalence.
\end{proposition}

\begin{proof}
Suppose that $P f$ is a strict weak equivalence.  Then it is also a weak
equivalence.  The maps $X \map P X$ and $Y \map P Y$ are weak
equivalences by Lemma \ref{lem:postnikov-we}, so $f$ is also.

Now suppose that $f$ is a weak equivalence.  By Proposition \ref{prop:tfae1},
we may assume that $f$ is a level map 
indexed by a cofinite directed set $I$ for which there
is a strictly increasing function $n: I \map \N$ such that 
$f_s: X_s \map Y_s$ is an $n(s)$-equivalence.

Consider the subsystem $X' = \{ P_{n(s)} X_s | s \in I \}$ of $X$
and the subsystem $Y' = \{ P_{n(s)} Y_s | s \in I \}$ of $Y$.  Note that
$X'$ and $Y'$ are cofinal in $X$ and $Y$.  Let $f'$ be the 
level map $X' \map Y'$ induced by $f$, so $f'$ is isomorphic to $f$.  
Since $X_s \map Y_s$ is an $n(s)$-equivalence, the map 
$P_{n(s)}X_s \map P_{n(s)}Y_s$ is a weak equivalence.  Hence
$f'$ is a levelwise weak equivalence, so $f$ is a strict weak equivalence.
\end{proof}

\begin{proposition} 
\label{prop:tfae3}
A map of pro-spaces $f: X \map Y$ is a weak equivalence if and only if
$\pi_0 f$ is an isomorphism of pro-sets, $\Pi_1 X \map f^*\Pi_1 Y$
is an isomorphism of pro-local systems on $X$, and for all $m$ and all
local systems $L$ on $Y$,
the map $H^m(Y; L) \map H^m(X; f^*L)$ is an isomorphism.
\end{proposition}

\begin{proof}
Let $f$ be a weak equivalence.  By Proposition \ref{prop:tfae1},
we may assume that 
$f$ is a level map indexed by a cofinite directed set $I$ for which there is
an increasing function $n: I \map \N$ such that $f_s$ is an
$n(s)$-equivalence.  

By the Whitehead theorem, $f_s$ induces a cohomology isomorphism
in dimensions less than $n(s)$ for any local system on $Y_s$.  Hence
$f$ induces an isomorphism
$H^m(Y; L) \map H^m(X; f^*L)$ in the colimit
for every $m$.   

Now suppose that $f$ satisfies the conditions of the proposition.
Factor $f$ as 

$$\xymatrix{X \ar[r]^i & Y' \ar[r]^p & Y,}$$
where $i$ is
a cofibration and $p$ is a strictly trivial fibration.
Since $p$ induces cohomology
isomorphisms by the first part of the proof,
the map $i$ still satisfies the hypotheses of
the proposition.  Therefore, 
we may assume that $f$ is a level map that is a level cofibration.

Note that $M = (X \times \Delta^1) \amalg_{X} Y$
is weakly equivalent to $Y$ since
$M$ is constructed levelwise and $M$ is levelwise weakly equivalent to $Y$.

We prove the proposition
by showing that for every strongly fibrant pro-space $Z$,
the map $\Map(M, Z) \map \Map(X, Z)$ is a weak equivalence.
A retract argument then shows that the map
$\Map(M, Z) \map \Map(X, Z)$ is a weak equivalence for all fibrant
pro-spaces $Z$.

Assume that $Z$ is an arbitrary strongly fibrant pro-space.  Note
that each $Z_s$ is a fibrant simplicial set with only finitely many
non-zero homotopy groups.  

Recall that $\Map(X,Z) = \lim_s \Map(X,Z_s)$.  Also recall that for every
$t$, the map $Z_t \map \lim_{s \leq t} Z_s$
is a fibration since $Z$ is fibrant.
Therefore, the map

$$\Map(X,Z_t) \map \lim_{s \leq t} \Map(X,Z_s) = \Map(X, \lim_{s \leq t} Z_s)$$
is a fibration.  It follows that
$\Map(X,Z)$ is weakly equivalent to the homotopy limit
$\holim_s \Map(X,Z_s)$.  Similarly,
$\Map(M,Z)$ is weakly equivalent to the homotopy limit
$\holim_s \Map(M,Z_s)$.

Since homotopy limits are invariant under levelwise weak equivalence,
we only need show that $\Map(M,Z_s) \map \Map(X,Z_s)$
is a weak equivalence of simplicial sets for each $s$.  Therefore, 
we may assume that 
$Z$ is a fibrant simplicial set with only finitely many non-zero homotopy
groups.

By adjointness, to show that $\Map (M, Z) \map \Map (X, Z)$
is a weak equivalence, it suffices to find lifts in the diagrams of
pro-spaces

$$\xymatrix{
X \ar[d]\ar[r] & Z^{\Delta^k} \ar[d] \\
M \ar[r]\ar@{-->}[ur] & Z^{\partial \Delta^k}.   }$$
For simplicity, rewrite the fibration
$Z^{\Delta^k} \map Z^{\partial \Delta^k}$ as $p: E \map B$.
Our goal is to find an $s$ and a commuting diagram of simplicial sets

$$\xymatrix{
X_s \ar[r]\ar[d] & E \ar[d] \\
M_s \ar[r]\ar[ur] & B.  }$$

Note that $E$ and $B$ are fibrant simplicial sets, and $p$ is a fibration
of simplicial sets.
Use Moore-Postnikov systems \cite[8.9]{JPM} to factor $E \map B$ as

$$\xymatrix{
E \ar[r] & \cdots \ar[r] & E_n \ar[r]^{p_n} & E_{n-1} \ar[r] &
\cdots \ar[r] & E_2 \ar[r]^{p_2} & E_1 \ar[r]^{p_1} & E_0 \ar[r]^{p_0} & B,}$$
where $E = \lim_n E_n$, each $p_n$ is a fibration, and the fibers of
$p_n$ (which may vary up to homotopy because $E_{n-1}$ may not be connected)
are of the form $K(\pi_n F, n)$ for some fiber $F$ of $p$.

The fibers of $p$ are of the form $\Omega^k Z$,
so they have
nonzero homotopy groups only in finitely many dimensions.  Therefore, there
exists $N$ such that $p_n$ is a weak equivalence for $n \geq N$.

We inductively construct partial liftings

$$\xymatrix{
X_s \ar[r]\ar[d] & E_n \ar[d] \\
M_s \ar[r]\ar[ur] & E_{n-1}.  }$$
Since $E = \lim_n E_n$, these lifts assemble to give us the desired lifting.

Choose $t$ such that the original square of pro-spaces
is represented by the square of simplicial sets

$$\xymatrix{
X_t \ar[r]\ar[d] & E \ar[d] \\
M_t \ar[r]       & B.   }$$
By the arguments of Lemma \ref{lem:technical-lift}
and Proposition \ref{prop:lift-trcofib-strong-fib},
there exists an $s_0 \geq t$ and a commutative diagram

$$\xymatrix{
X_{s_0} \ar[r]\ar[d] & E_0 \ar[d]^{p_0} \\
M_{s_0} \ar[r]\ar[ur]   & B.   }$$
These arguments apply because $p_0$ is a co-$1$-fibration and because
$X \map M$ induces an isomorphism in homotopy of dimension $0$ and $1$.

For the moment, assume that $k \geq 1$.  Since the fibers of $p$ are
of the form $\Omega^k Z$, $\pi_1 F$ is abelian for every fiber
$F$ of $p$.

By obstruction theory (for example, \cite{GWW}),
the only obstruction to finding a lift

$$\xymatrix{
X_{s_0} \ar[r]\ar[d] & E_1 \ar[d]^{p_1} \\
M_{s_0} \ar[r]\ar@{-->}[ur]   & E_0   }$$
is a class $\alpha$ in $H^2(M_{s_0}, X_{s_0}; L)$, where $L$ is the
local system given by the first homotopy groups of the fibers of $p_1$.  
Note that $L$ makes sense because these homotopy groups are abelian.

But $H^2(M,X; L) = 0$ since $X \map M$ induces a cohomology isomorphism.
Therefore, there exists an $s_1 \geq s_0$ such that the image of
$\alpha$ in $H^2(M_{s_1}, X_{s_1}; L)$ is $0$.  Hence, there is a lift
in the above diagram when $s_0$ is replaced by $s_1$.

The same obstruction theory argument applies inductively to give liftings
$$\xymatrix{
X_{s_n} \ar[r]\ar[d] & E_n \ar[d]^{p_n} \\
M_{s_n} \ar[r]\ar[ur]   & E_{n-1}   }$$
for each $ 2 \leq n \leq N$.  Now $L$ is the local system given by $n$th
homotopy groups of the fibers of $p_n$.  Again, $L$ makes sense because
$\pi_1 F = 0$ acts trivially on $\pi_n F$ for every fiber $F$ of $p_n$.  

Let $s = s_N$.  Recall that $p_n$ is a trivial fibration for $n \geq N$.
Therefore, lifts exist inductively in the diagrams
$$\xymatrix{
X_s \ar[r]\ar[d] & E_n \ar[d]^{p_n} \\
M_s \ar[r]\ar[ur]   & E_{n-1}   }$$
for $n > N$.  Hence, the desired lifting exists when $k \geq 1$.

Now consider $k = 0$.  The argument given for $k \geq 1$ does not work.
The trouble is that we cannot lift over $p_1$
with obstruction theory because the first homotopy groups of the fiber
are not necessarily abelian.

When $k = 0$, the map $p$ is just the map $Z \map *$, so we need to find
an $s$ and a factorization of $X_s \map Z$ through $M_s$.
Note that such factorizations are the same as factorizations up to 
homotopy of $X_s \map Z$ through $Y_s$.

Artin and Mazur \cite[Section 4]{AM} constructed such factorizations when
$Z$ is connected.  Their argument works even
when $Z$ is not connected provided that $\pi_0 X \cong \pi_0 Y$.  Here we
use the fact that $[X, Z]_{\pros} = \Hom_{\pro \Ho(\SSet)}(X, Z)$ by 
Lemma \ref{lem:into-constant}.

This proves the result.
\end{proof}

\section{Non-Cofibrantly Generated Model Structures} 
\label{sctn:cofibgen}

We prove in this section that the model structure of Section
\ref{sctn:structure} is not cofibrantly generated.  
The same argument shows that the strict model structure \cite{EH} is also not 
cofibrantly generated.  See Section \ref{sctn:strict} for a description
of the strict structure.  We start with a general lemma about
cofibrantly generated model structures.

\begin{lemma} 
\label{lem:cofibgen}
Suppose that a model structure on a category $\C$ is cofibrantly generated
with a set of generating cofibrations $I$.  Let $T$ be the set of targets
of maps in $I$, and let $X$ be any cofibrant object of $\C$
not isomorphic to the initial object.  Then there exists some
$Y$ in $T$ not isomorphic to the initial object with a map $Y \map X$
in $\C$.
\end{lemma}

\begin{proof}
Let $X$ be a cofibrant object of $\C$.
Then $X$ is a retract of another object
$X'$, where $X'$ is a transfinite composition of pushouts of maps in $I$
\cite[14.2.12]{PH}.
Since there is a map from $X'$ to $X$, it suffices to find a map from 
some object of $T$ to $X'$.  Since $X$ is not the initial object,
$X'$ is also not the initial object.  Hence $X'$ is a non-trivial
transfinite composition of pushouts of maps in $I$.  Let $Z \map Y$ be
a map in $I$ occurring in the construction of $X'$.  Then there is a
map from $Y$ to $X'$.
\end{proof}

The next proposition gives a construction of specific pro-sets with 
special properties.

\begin{proposition} 
\label{prop:large-objects}
Let $F$ be a small family of pro-sets ({\em i.e.}, a set of pro-sets) 
not containing the empty pro-set.
Then there exists a pro-set $X$ such that for every $Y$ in $F$, there
are no maps $Y \map X$ of pro-sets.
\end{proposition}

\begin{proof}
Choose an infinite
cardinal $\kappa$ larger than the size of any of the sets occurring
in any of the objects of $F$.  Let $S$ be a set of size $\kappa$.

Define a pro-set $X$ as follows.  
Consider the collection of all subsets $U$ of $S$ whose complements
$U^c$ are strictly smaller than $S$.  Note that this implies that the
size of $U$ is $\kappa$, but the converse is not true.
These subsets form a pro-set, where the structure maps are
inclusions.  This system is cofiltered because $(U \cap V)^c =
U^c \cup V^c$ is strictly smaller than $S$ when $U^c$ and $V^c$ are.

Let $Y$ be an object of $F$.  Suppose that there is a map $f: Y \map X$ of
pro-sets.  Then there exists a $t$ and a map $f_{t,S}: Y_t \map S$ 
representing $f$.  Let $A$ be the image of $f_{t,S}$, so 
$A$ is strictly smaller than $S$ since $Y_t$ is strictly smaller than
$S$.  Consider the set $S - A$, which occurs as an object in the system
$X$.  Since $f$
is a map of pro-sets, there exists a $u \geq t$ such that
the composition $Y_u \map Y_t \map S$ factors through $S - A$.
Since $Y_t$ and $S - A$ have disjoint images in $S$, 
this is only possible if $Y_u$ is the empty set.  However, $Y_u$ cannot
be the empty set because $Y$ is not the empty pro-set.  By contradiction,
the map $f$ cannot exist.
\end{proof}

\begin{corollary} 
\label{cor:uncofibgen}
There are no cofibrantly generated model structures on pro-spaces for
which every object is cofibrant.
\end{corollary}

\begin{proof}
We argue by contradiction.
Suppose that there exists a cofibrantly generated model structure for which
every object is cofibrant.  
Let $I$ be the set of generating cofibrations, and let $T$ be the
set of targets of maps in $I$.  Apply $\pi_0$ to $T$ to obtain a small family
of pro-sets $F$.

Let $X$ be the pro-set constructed in Proposition \ref{prop:large-objects}.
We can
think of $X$ as a pro-space by identifying a set with a simplical set of 
dimension zero.  By Lemma \ref{lem:cofibgen}, there exists a non-empty $Y$ in
$T$ and a map $Y \map X$.  This induces a map $\pi_0 Y \map \pi_0 X = X$.
However, such a map cannot exist by Proposition \ref{prop:large-objects}
because $\pi_0 Y$ belongs to $F$.
\end{proof}

This corollary applies in particular to the model structure of 
Section \ref{sctn:structure} and to the strict model structure.


\begin{thebibliography}{99}

\bibitem{SGA}
M. Artin, A. Grothendieck, and J. L. Verdier, {\em Theorie des topos
et cohomologie \'etale des schemas}, Lecture Notes in
Mathematics, vol. 269, Springer Verlag, 1972.

\bibitem{AM}
M. Artin and B. Mazur, {\em Etale homotopy}, 
Lecture Notes in Mathematics, vol. 100, Springer Verlag, 1969.

\bibitem{BF}
A. K. Bousfield and E. M. Friedlander, {\em 
Homotopy theory of $\Gamma$-spaces,
spectra, and bisimplicial sets}, Geometric Applications of
Homotopy Theory, vol. II (Proc. Conf., Evanston, IL, 1977),
Lecture Notes in Mathematics, vol. 658, Springer Verlag, 1978, pp. 80--130.

\bibitem{BK}
A. K. Bousfield and D. Kan, {\em Homotopy limits, completions, and
localizations}, Lecture Notes in Mathematics, vol. 304, Springer
Verlag, 1972.

\bibitem{etK}
W. G. Dwyer and E. M. Friedlander, {\em Algebraic and \'etale $K$-theory},
Trans. Amer. Math. Soc. {\bf 292} (1985), 247--280.

\bibitem{EH}
D. A. Edwards and H. M. Hastings, {\em Cech and Steenrod homotopy
theories with applications to geometric topology}, Lecture Notes in
Mathematics, vol. 542, Springer Verlag, 1976.

\bibitem{etK-1}
E. M. Friedlander, {\em Etale $K$-theory I:
Connections with \'etale cohomology and algebraic vector bundles},
Invent. Math. {\bf 60} (1980), 105--134.

\bibitem{EF}
\bysame, {\em Etale homotopy of simplicial schemes}, Annals of
Mathematics Studies, vol. 104, Princeton University Press, 1982.

\bibitem{JG}
J. Grossman, {\em A homotopy theory of pro-spaces}, Trans. Amer. Math.
Soc. {\bf 201} (1975), 161--176.

\bibitem{PH}
P. Hirschhorn, {\em Localization of Model Categories}, preprint.

\bibitem{MAM}
M. A. Mandell, {\em $E_\infty$ algebras and $p$-adic homotopy theory},
to appear in Topology.

\bibitem{JPM}
J. P. May, {\em Simplicial objects in algebraic topology}, Van Nostrand
Mathematical Studies, vol. 11, Van Nostrand, 1967.

\bibitem{FM}
F. Morel, {\em Ensembles profinis simpliciaux et interpr\'etation 
g\'eom\'etrique du foncteur $T$}, Bull. Soc. Math. France {\bf 124}
(1996), 347--373.

\bibitem{DQ}
D. G. Quillen, {\em Homotopical algebra},
Lecture Notes in Mathematics, vol. 43, Springer Verlag, 1967.

\bibitem{DS}
D. Sullivan, {\em Genetics of homotopy theory and the Adams conjecture},
Ann. of Math. {\bf 100} (1974), 1--79.

\bibitem{GWW}
G. W. Whitehead, {\em Elements of homotopy theory}, Graduate
Texts in Mathematics, vol. 61, Springer Verlag, 1978.

\end{thebibliography}
\end{document}